# The Computation of the Mean First Passage Times for Markov Chains

Jeffrey J. Hunter[1]

*Department Mathematical Sciences*
*School of Engineering, Computer and Mathematical Sciences, Auckland University of Technology,*
*Private Bag 92006, Auckland 1142, New Zealand*

**Abstract**

A survey of a variety of computational procedures for finding the mean first passage times in Markov chains is presented. The author recently developed a new accurate computational technique, an Extended *GTH* Procedure, Hunter (Special Matrices, 2016) similar to that developed by Kohlas (Zeit. fur Oper. Res., 1986). In addition, the author recently developed a variety of new perturbation techniques for finding key properties of Markov chains including finding the mean first passage times, Hunter (Linear Algebra and its Applications, 2016). These recently developed procedures are compared with other procedures including the standard matrix inversion technique using the fundamental matrix (Kemeny and Snell, 1960), some simple generalized matrix inverse techniques developed by Hunter (Asia Pacific J. Oper. Res., 2007), and some modifications to the *FUND* technique of Heyman (SIAM J Matrix Anal. and Appl., 1995). *MATLAB* is used to compute errors and estimate computation times when the techniques are used on some test problems that have been used in the literature together with some large sparse state-space cases. For accuracy a preference for the procedure of the author is exhibited for the test problems. However it appears that the procedure, as presented, requires longer computational times.

## 1. Introduction

In Markov chain (*MC*) theory mean first passage times (*MFPT*s) provide significant information regarding the short term behaviour of the *MC*. A review of *MFPT*s, together with details regarding stationary distributions and the group inverse of the Markovian

---

kernel, is given in [18]. We refer the reader to this aforementioned article as it provides the relevant background to this paper and enables us to avoid repetition of the material. In Hunter [18], which focuses on computational techniques for the key properties of irreducible *MC*s using perturbation techniques, we commented that in a sequel paper we would consider a variety of other techniques to get a better impression as to whether perturbation procedures may in fact prove to be suitable alternatives. We address these issues in this paper.

We firstly set the scene by reintroducing the notation that was used in [18].

Let $\{X_n, n \geq 0\}$ be a finite *MC* with state-space $S = \{1, 2, ..., m\}$ and transition matrix $P = [p_{ij}]$, where $p_{ij} = P\{X_n = j \mid X_{n-1} = i\}$ for all $i, j \in S$.

The stationary distribution $\{\pi_j\}$, $(1 \leq j \leq m)$, exists and is unique for all irreducible *MC*s, that $\pi_j > 0$ for all $j$, and satisfies the equations (the *stationary equations*)

$$\pi_j = \sum_{i=1}^{m} \pi_i p_{ij} \text{ with } \sum_{i=1}^{m} \pi_j = 1. \tag{1.1}$$

If $\pi^T \equiv (\pi_1, \pi_2, ..., \pi_m)$, the *stationary probability vector*, and $e$ is a column vector of 1's, the stationary equations (1.1) can be expressed as

$$\pi^T(I - P) = \mathbf{0}^T, \text{ with } \pi^T e = 1. \tag{1.2}$$

Let $T_{ij} = \min[n \geq 1, X_n = j \mid X_0 = i]$ be the first passage time from state $i$ to state $j$ (first return when $i = j$) and define $m_{ij} = E[T_{ij} \mid X_0 = i]$ as the *MFPT* from state $i$ to state $j$ (or mean recurrence time of state $i$ when $i = j$). For finite irreducible *MC*s all the $m_{ij}$ are well defined and finite. Let $M = [m_{ij}]$ be the *MFPT* matrix. Let $\delta_{ij} = 1$, when $i = j$ and 0, when $i \neq j$. Let $M_d = [\delta_{ij} m_{ij}]$ be the diagonal matrix formed from the diagonal elements of $M$, and $E = [1]$ (i.e. all the elements are unity). Let $\Pi = e\pi^T$.

It is well known ([19]) that, for $1 \leq i, j \leq m$,

$$m_{ij} = 1 + \sum_{k \neq j} p_{ik} m_{kj}. \tag{1.3}$$

In particular, for all $j \in S$, the mean recurrence time of state $j$ is given by

$$m_{jj} = 1/\pi_j. \tag{1.4}$$

From (1.3) and (1.4) it follows that $M$ satisfies the matrix equation

$$(I - P)M = E - PM_d, \text{ with } M_d = (\Pi_d)^{-1}. \tag{1.5}$$

Note that the expression (1.5) typically involves knowledge of $P_d$, i.e. the stationary probabilities. In this paper, as we are not focussing on the computation of stationary distributions, when we require such terms when they are not explicitly derived in carrying out the computations for $M$, we typically use the *GTH* algorithm of Grassman, Taksar and Heyman [3], (or equivalently the State Reduction procedure of Sheskin [23]), as these are known to give accurate results with no subtractions being involved. There are however other alternative procedures that could be used, for example the *eig* procedure of *MATLAB*.

We provide twelve procedures for solving, in effect, equations (1.3) or (1.5), for the *MFPT*s. In Section 2 we give some direct procedures (Procs 1 and 2) based upon utilising matrix inverses. In Section 3 we summarise the six perturbation procedures (Procs 3 to 8), given in Hunter [18]. In Section 4 we outline the extended *GTH* procedure (*EGTH*) of Hunter [17] based upon Kohlas [20] (Proc 9), while in Section 5 we outline, modify and simplify the *FUND* procedure of Heyman [7], putting it in a generalized matrix inverse framework, that enables us to find the *MFPTs* without directly computing the fundamental matrix, to yield Procs 10, 11 and 12.

In Section 6 we describe a set of test problems, used initially by Harrod and Plemmons [5] in comparing different techniques for computing the stationary probabilities. These are augmented with two randomly generated large sparse state-space examples.

In Section 7 we use *MATLAB* computations to compare errors. We compare, in particular, the percentage of zero errors, the overall residual errors of our computations, in double precision, as well as the number additional accurate digits achieved with double precision over single precision. This leads to the conclusion that, typically, the *EGTH* Procedure of Hunter [17] gives us the most accurate results especially for the prescribed small state-space test problems.

In Section 8 we compare the computational times of the procedures using the "tic-toc" procedure of *MATLAB*, averaged over a number of test runs. The surprising observation is that the *EGTH* Procedure is by far the most time consuming technique but this is due to the staged implementation of the procedure as programmed with some computed *MFPTs* being omitted.

Section 9 gives some conclusions and summarises the results to provide some implementation guidelines and recommendations.

## 2. Computation of *MFPTs* using matrix inverses

If $A$ is an $m \times m$ matrix of real elements and $X$ is any $m \times m$ matrix that satisfies the condition $AXA = A$, then $X$ is said to be a one-condition generalised matrix inverse, a g-inverse, of $A$, and is often written as $A^-$. If $A$ is non-singular then $A^- = A^{-1}$.

All g-inverses of $I - P$ can be expressed in terms utilising matrix inverses, as pioneered by Hunter [11]. The general result is as follows:

**Theorem 2.1**: *Let $P$ be the transition matrix of a finite irreducible Markov chain with m states and stationary probability vector $\boldsymbol{\pi}^T = (\pi_1, \pi_2, \ldots, \pi_m)$. Let $\boldsymbol{e}^T = (1, 1, \ldots, 1)$ and $\boldsymbol{t}$ and $\boldsymbol{u}$ be any vectors.*

(a)   $I - P + \boldsymbol{tu}^T$ *is non-singular if and only if $\boldsymbol{\pi}^T \boldsymbol{t} \neq 0$ and $\boldsymbol{u}^T \boldsymbol{e} \neq 0$.*

(b)   *If $\boldsymbol{\pi}^T \boldsymbol{t} \neq 0$ and $\boldsymbol{u}^T \boldsymbol{e} \neq 0$ then $[I - P + \boldsymbol{tu}^T]^{-1}$ is a one-condition g-inverse of $I - P$.*

(c)   *All one-condition g-inverses of $I - P$ can be expressed as*
$$A^- = [I - P + \boldsymbol{tu}^T]^{-1} + \boldsymbol{ef}^T + \boldsymbol{g\pi}^T \text{ for arbitrary vectors } \boldsymbol{f} \text{ and } \boldsymbol{g}.$$
□

Well-known special g-inverses of $I - P$ are Kemeny and Snell's *fundamental matrix* $Z = [I - P + \Pi]^{-1}$ where $\Pi = \boldsymbol{e\pi}^T$, introduced in [19], (and initially shown to be a g-inverse of $I - P$ by Hunter [10]) and Meyer's *group inverse of $I - P$* given by $A^{\#} = Z - \Pi$ ([21]). ($A^{\#}$ is more restrictive than a simple g-inverse in that it is the unique g-inverse that satisfies $(I - P)A^{\#} = A^{\#}(I - P) = I - e\pi^T$, $A^{\#}e = 0$ and $\pi^T A^{\#} = 0^T$ where $A = I - P$.)

One-condition g-inverses of $I - P$ are introduced as they are typically used to solve systems of linear equations involving $I - P$ (as in (1.2) and (1.5)).

Solving the equations given by (1.5), using Theorem 2.1, yields the following general results for finding the *MFPTs* of *MC's*, (see [11], [12] for (a) and [16] for (b) and (c)):

**Theorem 2.2**:
*(a) If G is any g-inverse of $I - P$, then the MFPT matrix M, is given as*
$$M = [G\Pi - E(G\Pi)_d + I - G + EG_d]D, \quad (2.1)$$
where $D = (\Pi_d)^{-1} = \left[(e\pi^T)_d\right]^{-1}$.
*(b) If $H \equiv G(I - P)$ then H is a g-inverse of $I - P$ with $He = 0$ and*
$$M = [I - H + EH_d]D. \quad (2.2)$$
*(c) $Ge = ge$ for some g if and only if*
$$M = [I - G + EG_d]D. \quad (2.3)$$
□

Special cases of (2.3) for $M$ are $G = Z$ and $G = A^{\#}$.

Theorem 2.2 above leads to the following two procedures.

**Proc 1: (Standard method)**
*Given an irreducible P*
*(i) Compute the stationary probability vector $p^T$.*
*(ii) Compute the fundamental matrix $Z = [I - P + e\pi^T]^{-1}$.*
*(iii) Compute $M = [I - Z + EZ_d][(e\pi)_d]^{-1}$.*

This is the original procedure developed by Kemeny and Snell [19] and has been universally used in the past. As identified above, prior to computing Z, the stationary probability vector $\pi^T$ is required. We use the *GTH* algorithm to compute $\pi^T$. (See Section 4 for details.)

Hunter [15] established a number of results regarding expressions for the *MFPTs* using a range of simple g-inverses of the form given in Theorem 2.1, (typically with *f* and *g* taken as zero vectors.) The simplest result is given as follows.

**Proc 2: (Simple method)**
*Given an irreducible P*
*(i) Compute the g-inverse $G = \left[I - P + ee_b^T\right]^{-1}$ where $e_b^T$ is a vector with 1 in the b-th position and 0 elsewhere.*
*(ii) Compute, $\pi^T = Ge_b^T$ so that if $G = [g_{ij}]$ then $\pi_j = g_{bj}$, $j = 1, 2, ..., m$.*
*(iii) Compute $M = [I - G + EG_d][(e\pi^T)_d]^{-1}$.*

Thus following one matrix inversion, one can find the stationary probabilities (requiring only the $b^{th}$ row) and the mean first passage times. The choice of $b$ is arbitrary. We take $b = 1$ in our test examples (in Section 6).

Since the above two procedures both require the evaluation of a matrix inverse we do not expect them to perform well in examples when we have either a large number of states or ill-conditioned matrices.

## 3. Computation of *MFPTs* using perturbation procedures

The general idea behind the perturbation procedures, which are considered in detail in Hunter [18], is the following. Start with a simple transition matrix $P_0$ with known or easily computed stationary probability vector $\pi_0^T$, mean first passage time matrix $M_0$ and a simple g-inverse $G_0$ (or easily computed fundamental matrix $Z_0$ or group inverse $A_0^\#$.) Update the transition matrices $P_{i-1}$ by sequentially replacing the *i-th* row of $P_{i-1}$ with $p_i^T = e_i^T P$, the *i-th* row of the given transition matrix $P$ to yield $P_i$, ($i = 1, 2, …, m$) ending up with $P_m = P$. The simplest structure to start with is the irreducible transition matrix $P_0 = ee^T/m$ as this ensures that each subsequent updated transition matrix is also irreducible. As in [18], if $P = \sum_{i=1}^{m} e_i p_i^T$ then $P_i = P_{i-1} + e_i b_i^T$ with $b_i^T = p_i^T - e^T/m$, for $i = 1, 2, …, m$. We update $\pi_{i-1}^T$, $M_{i-1}$, and $G_{i-1}$ ( or $Z_{i-1}$, $A_{i-1}^\#$) to $\pi_i^T$, $M_i$, and $G_i$ ( or $Z_i$, $A_i^\#$) finishing with $\pi^T = \pi_m^T$, $M = M_m$ and $G = G_m$, (or $Z = Z_m$, $A^\# = A_m^\#$). With $P_0$ as above, $\pi_0^T = e^T/m$, $Z_0 = I$, $A_0^\# = I - ee^T/m$ and $M_0 = mee^T$.

The successive updates effectively make use of the Sherman-Morrison [22] formula for computing matrix inverses. The details are given in [18].

The first perturbation procedure is an extension to the procedure of Hunter [14] where an updated one-condition g-inverse that is used to find successive stationary probability vectors is utilised to compute the *MFPT* matrix. Let $G_i = [I - P_i + t_i u_i^T]^{-1}$. We update the g-inverse $G_{i-1}$ to $G_i$ successively with $t_0 = e$, $u_0^T = e^T/m$ (i.e. $G_0 = [I - P_0 + t_0 u_0^T]^{-1} = I$) and $t_i = e_i$, $u_i^T = u_{i-1}^T + b_i^T$ ($i = 1, …, m$). We use Theorem 2.2(b) as this eliminates the requirement to find the group inverse but utilises the structure of $H$, a particular g-inverse of $I - P$, to find $M$. This is Algorithm 1 in [18].

**Proc 3: (G-inverse update – Pert AL1)**

    *(i)*    Let $G_0 = I$, $u_0^T = e^T/m$.

    *(ii)*   For $i = 1, 2, ..., m$, let $p_i^T = e_i^T P$, $u_i^T = u_{i-1}^T + p_i^T - e^T/m$,

             $G_i = G_{i-1} + G_{i-1}(e_{i-1} - e_i)(u_{i-1}^T G_{i-1} / u_{i-1}^T G_{i-1} e_i)$.

    *(iii)*  At $i = m$, let $G_m = G$ and $\pi^T = \pi_m^T = \dfrac{u_m^T G_m}{u_m^T G_m e}$.

    *(iv)*  Compute $H = G(I - e\pi^T)$.

(vi) Compute $M = [I - H + E(diag(H))]D$ where
$E = [1]$ and $D = inv[diag(e\pi^T)]$.

For Proc 4 we consider a modification of Algorithm 2 given in [18]. We consider successive row perturbations of the group inverse $A_i^\#$ leading to $A^\#$. However we do not need to compute the group inverse $A_i^\#$ at each stage but rather $R_i$ where $A_i^\# = R_i + ey_i^T$. We start with $A_0^\# = R_0 = I - ee^T/m$ leading to $R_m$ which is a g-inverse of $I - P$ with the property that $R_m e = 0$. We utilise Theorem 2.2(c) to compute $M$.

**Proc 4: (Modified group inverse update – Pert AL2)**
Start with $P$.
(i) Set $R_0 = I - ee^T/m$.
(ii) For $i = 1, 2, \ldots, m$, let $p_i^T = e_i^T P$, $b_i^T = p_i^T - e^T/m$,
$$R_i = R_{i-1} + \frac{1}{1 - b_i^T R_{i-1} e_i} R_{i-1} e_i b_i^T R_{i-1}.$$
(iii) Compute $\pi^T = e_1^T - e_1^T (I - P) R_m$.
(iv) Compute $M = [I - R_m + E(diag(R_m))]D$, where $E = [1]$ and $D = inv[diag(e\pi^T)]$.

In [18] it is shown that not all the calculations are required. At the $i$-th recursion, leading to $R_i$, the only terms that are updated are in the first $i$ rows with the rows numbered $i+1$, $i+2, \ldots, m$ remaining unchanged.

For Proc 5, rather than focus directly on row operations to lead to the group inverse and thence $M$, one can develop a matrix procedure that leads jointly to the matrix $\Pi = e\pi^T$ and the group inverse $A^\#$.

Under a perturbation $E$ when $\pi^T$ leads to $\bar{\pi}^T = (\pi^T(I - EA^\#)^{-1})$, if $\Pi = e\pi^T$ and $\bar{\Pi} = e\bar{\pi}^T$ then $\bar{\Pi} = \Pi(I - EA^\#)^{-1}$.

Under the perturbation $E = e_i b_i^T$ to the $i$-th row with $b_i^T e = 0$, yields,
$$\bar{\Pi} = \Pi\left[I + \frac{1}{1 - b_i^T A^\# e_i} e_i b_i^T A^\#\right] \text{ and } \bar{A}^\# = (I - \bar{\Pi})A^\#\left(I + \frac{1}{1 - b_i^T A^\# e_i} e_i b_i^T A^\#\right).$$

This leads to the following procedure. (For more details see Algorithm 3 in [18]).

**Proc 5: (Group inverse by matrix updating – Pert AL3)**
(i) Let $P_0 = ee^T/m$, implying $\Pi_0 = ee^T/m$, $A_0^\# = I - ee^T/m$.
(ii) For $i = 1, 2, \ldots, m$, let $p_i^T = e_i^T P$, $b_i^T = p_i^T - e^T/m$,
$$S_i = I + \frac{1}{1 - b_i^T A_{i-1}^\# e_i} e_i b_i^T A_{i-1}^\#, \quad \Pi_i = \Pi_{i-1} S_i, \quad A_i^\# = (I - \Pi_i)A_{i-1}^\# S_i.$$
(iii) At $i = m$, let $S = S_m$ then $\Pi = \Pi_{m-1} S$, $A^\# = (I - \Pi)A_{m-1}^\# S$.

*(iv)* Compute $M = [I - A^\# + EA_d^\#]D$, where $E = [1]$ and $D = (\Pi_d)^{-1}$.

Procs 6, 7 and 8 to follow are three interrelated procedures, each with different starting conditions, based on updating simple g-inverses of $I - P_0$ that lead to simple computations for the stationary probabilities and the *MFPT* matrix.

From Theorem 2.2(c), if we choose a g-inverse $G$ of $I - P$ with the property that $Ge = ge$, by taking $G$ of the form $G = [I - P + e\beta^T]^{-1}$ then $\pi^T = \beta^T G$. Further we have a simple form of the *MFPT* matrix $M$ given by eqn. (2.3). (While it is easy to find an expression for the group inverse of $I - P$ as $A^\# = (I - e\pi^T)G$ we don't actually require that step to find expressions for $M$.)

In Hunter [15] we explored the properties of some g-inverses of this form. For the three procedures to follow we use, successively, the special forms, $G_e \equiv \left[I - P + \dfrac{ee^T}{m}\right]^{-1}$, $G_{e1} \equiv [I - P + ee_1^T]^{-1}$ and $G_{ee} \equiv [I - P + ee^T]^{-1}$, and the Sherman-Morrison [22] matrix inversion formula. The starting conditions for each procedure are different and, although we carry out similar recursions, we have different expressions for the stationary probability vector $p^T$ but identical calculation procedures for the *MFPT* matrices. In each procedure $K_i = [I - P_i + e\beta^T]^{-1}$ with $K_0$ as specified leading to $K_m$ as the required matrix inverse. See [18] for full details.

This leads to three further algorithms – Algorithms 4A, 4B, and 4C in [18]. They are all variants of the generic recursion given in Proc 6, with identical steps *(ii)* and *(iv)* but different initial conditions *(i)* and final step at the *m*-th iteration *(iii)*.

**Proc 6: (Update using $G_e$ - Pert AL4A)**
*(i)* Start with $K_0 = I$.
*(ii)* For $i = 1, 2, ..., m$, let $p_i^T = e_i^T P$, $b_i^T = p_i^T - e^T/m$,

$$K_i = K_{i-1}(I + C_i), \text{ where } k_i = 1 - b_i^T K_{i-1} e_i \text{ and } C_i = \dfrac{1}{k_i} e_i b_i^T K_{i-1}.$$

*(iii)* At $i = m$, let $K = K_m$ and then compute $\pi^T = \dfrac{1}{m} e^T K$.

*(iv)* Compute $M = [I - K + EK_d]D$, where $E = [1]$ and $D = \left[(e\pi^T)_d\right]^{-1}$.

**Proc 7: (Update using $G_{e1}$ - Pert AL4B)**
*(i)* Start with $K_0 = I + e\left(\dfrac{e^T}{m} - e_1^T\right)$.
*(ii)* Carry out Step *(ii)* of Proc 6.
*(iii)* At $i = m$, let $K = K_m$ and then compute $\pi^T = e_1^T K$.
*(iv)* Carry out Step *(iv)* of Proc 6, to compute $M$.

**Proc 8: (Update using $G_{ee}$ – Pert AL4C)**

*(i)* Start with $K_0 = I - \left(\dfrac{m-1}{m^2}\right)ee^T$.

*(ii)* Carry out Step *(ii)* of Proc 6.

*(iii)* At $i = m$, let $K = K_m$ and then compute $\pi^T = e^T K$.

*(iv)* Carry out step *(iv)* of Proc 6, to compute $M$.

## 4. Computation of *MFPTs* using Hunter Extended *GTH* (*EGTH*) procedure

The details of this *EGTH* procedure are given in Hunter [17]. We make use of the *GTH* procedure of Grassman, Taksar and Heyman [3] (or the equivalent state reduction procedure by Sheskin [23]) for finding the stationary probability vector *p*. We provide some details that serves to introduce some additional notation.

Start with the given transition matrix $P \,(= P^{(m)})$ of the irreducible *MC* $\{X_k^{(m)}, k \geq 0\}$ with state-space $S = \{1,2,\ldots, m\} \equiv S_m$. The general idea is to reduce the state-space, one state at a time successively removing states $m-1, m-2, \ldots$ until we are left with a single state 1. Once state 1 is reached the state-space is expanded one state at a time i.e. inserting states $2, \ldots$, successively to finally insert state *m*.

Suppose we reach the stage where we have $n$ states $S_n = \{1, 2, \ldots n\}$ with *MC* $\{X_k^{(n)}, k \geq 0\}$ and transition matrix $P^{(n)}$, then it is easily shown during the state reduction process that the elements of $P^{(n-1)} = \left[p_{ij}^{(n-1)}\right]$ are related to the earlier elements of $P^{(n)}$ as

$$p_{ij}^{(n-1)} = p_{ij}^{(n)} + \frac{p_{in}^{(n)} p_{nj}^{(n)}}{S(n)}, \quad 1 \leq i \leq n-1,\ 1 \leq j \leq n-1, \tag{4.1}$$

where $S(n) = 1 - p_{nn}^{(n)} = \sum_{j=1}^{n-1} p_{nj}^{(n)}$. Note that the transition probabilities of the reduced *MC* can all be obtained without carrying out any subtraction. The *MC* $\{X_k^{(n-1)}, k \geq 0\}$ on the reduced state-space, $S_{n-1}$ is the "censored" *MC* (see [2]), i.e. the *MC* restricted to the states of $S_{n-1}$. Further, the irreducibility of the reduced state-space *MC* is retained. One can derive relationships between the stationary distributions of the respective *MCs*, i.e. $\{\pi_i^{(n)}, i \in S_n\}$ for $\{X_k^{(n)}, k \geq 0\}$ on $S_n$. In particular, it can be shown,

$$\pi_i^{(n-1)} = \frac{\pi_i^{(n)}}{1 - \pi_n^{(n)}} = \frac{\pi_i^{(n)}}{\sum_{k=1}^{n-1} \pi_k^{(n)}}, \quad 1 \leq i \leq n-1.$$

Similarly, when we expand the state-space we can show that

$$\boldsymbol{\pi}^{(n)} \equiv \left(\pi_1^{(n)}, \ldots, \pi_n^{(n)}\right) = c_{n-1}\left(\pi_1^{(n-1)}, \ldots, \pi_{n-1}^{(n-1)}, \frac{\sum_{i=1}^{n-1} \pi_i^{(n-1)} p_{in}^{(n)}}{S(n)}\right),$$

where $c_{n-1}$ is determined from the fact that $\sum_{i=1}^{n} \pi_i^{(n)} = 1$.

From these results we have the following algorithm.

**GTH Procedure for computing the stationary probabilities of a *MC*:**

Let MC $\{X_k^{(m)}, k \geq 0\}$ be finite irreducible MC with state-space $S_m = \{1, 2, \ldots, m\}$ and transition matrix $P = P^{(m)} = [p_{ij}^{(m)}]$. Let $\{\pi_i^{(m)}\}$ be its stationary probabilities.

Step 1. Compute, successively for $n = m, m-1, \ldots, 3$,

$$p_{ij}^{(n-1)} = p_{ij}^{(n)} + \frac{p_{in}^{(n)} p_{nj}^{(n)}}{S(n)}, \quad 1 \leq i \leq n-1, \, 1 \leq j \leq n-1 \text{ where } S(n) = \sum_{j=1}^{n-1} p_{nj}^{(n)}.$$

Step 2. Set $r_1 = 1$ and compute successively for $n = 2, \ldots, m$, $r_n = \sum_{i=1}^{n-1} r_i p_{in}^{(n)} / S(n)$.

Step 3. Compute, for $i = 1, 2, \ldots, m$, $\pi_i^{(m)} = r_i / \sum_{j=1}^m r_j$. □

In extending this algorithm to find the *MFPTs*, Kohlas [20] showed that it is more natural to consider the process as a Markov renewal process (*MRP*), $\{(X_k^{(n)}, T_k^{(n)}), k \geq 0\}$, where $\{X_k^{(n)}, k \geq 0\}$, is the embedded *MC* when the state-space is $S_n = \{1, \ldots, n\}$ and $T_k^{(n)}$ is time that the process stays in the state before making the next transition. Let $\mu_i^{(n)} = E[T_k^{(n+1)} - T_k^{(n)} | X_k^{(n)} = i]$ be the expected holding time in state $i$ when the state-space is $S_n$. When the process is censored by eliminating state $n$ the mean holding time vector eliminates that state and reduces to a smaller $(n-1)$-dimension vector as

$$\boldsymbol{\mu}^{(n-1)T} = (\mu_1^{(n-1)}, \ldots, \mu_{n-1}^{(n-1)}) \text{ where } \mu_i^{(n-1)} = \mu_i^{(n)} + \frac{p_{in}^{(n)} \mu_n^{(n)}}{S(n)}, \quad 1 \leq i \leq n-1.$$

Under the *MC* setting, which we assume in this paper, initially $\mu_i^{(m)} = 1$ for all $i \in S_m$. In [17] it is shown how this influences the *MFPTs* showing, in particular, that

$$m_{ij} = \frac{\mu_i^{(i)} + \sum_{k=1, k \neq j}^{i-1} p_{ik}^{(i)} m_{kj}}{S(i)}, \quad 3 \leq i \leq m, 1 \leq j \leq i-1,$$

with $m_{21} = \frac{\mu_2^{(2)}}{S(2)}$, and $m_{ii} = \mu_i^{(i)} + \sum_{k=1,}^{i-1} p_{ik}^{(i)} m_{ki}$, $2 \leq i \leq m$, with $m_{11} = \mu_1^{(1)}$.

The expressions for $m_{ij}$, for $1 \leq i \leq m-1, i+1 \leq j \leq m$ are much more complicated. However, by focussing primarily on the terms $m_{i1}$ for $1 \leq i \leq m$, i.e. the first column of the matrix of *MFPTs*, we can produce a simple algorithmic procedure.

**Proc 9: (*EGTH* – Hunter Extended *GTH* Procedure)**
Let $\{X_k^{(m)}, k \geq 0\}$ be a finite irreducible *MC* with state-space $S_m = \{1, 2, \ldots, m\}$ and transition matrix $P \equiv P^{(m)} \equiv [p_{ij}^{(m)}]$.

Step 1(*i*): Carry out step 1 of the *GTH* Procedure, i.e.
Compute, successively for $n = m, m-1, \ldots, 3$,

$$p_{ij}^{(n-1)} = p_{ij}^{(n)} + \frac{p_{in}^{(n)} p_{nj}^{(n)}}{S(n)}, \quad 1 \leq i \leq n-1, \, 1 \leq j \leq n-1 \text{ where } S(n) = \sum_{j=1}^{n-1} p_{nj}^{(n)}. \quad (4.2)$$

Step 1(*ii*): Compute, successively for $n = m, m-1, \ldots, 3, 2$,

$$\mu_i^{(n-1)} = \mu_i^{(n)} + \frac{\mu_n^{(n)} p_{in}^{(n)}}{S(n)}, \quad 1 \leq i \leq n-1, \text{ where } \boldsymbol{\mu}^{(m)T} = (\mu_1^{(m)}, \ldots, \mu_m^{(m)}) = (1, \ldots, 1).$$

Step 1(*iii*): Compute the $m \times 1$ column vector $\boldsymbol{m}_m^{(1)(m)} = (m_{i1})$,

where $m_{11} = \mu_1^{(1)}$, $m_{21} = \dfrac{\mu_2^{(2)}}{S(2)}$, and for $i = 3, ,,,, m$, $m_{i1} = \dfrac{\mu_i^{(i)} + \sum_{k=2}^{i-1} p_{ik}^{(i)} m_{k1}}{S(i)}$.

Thus, starting with $P^{(m)} \equiv P^{(m)(1)}$, we can easily obtain the entries of the first column of the matrix $M$ i.e. $\boldsymbol{m}_m^{(1)(m)}$, where $M = [m_{ij}] = (\boldsymbol{m}_m^{(1)(m)}, \boldsymbol{m}_m^{(2)(m)},...,\boldsymbol{m}_m^{(m)(m)})$.

The procedure that follows to find the other *MFPTs* is to permute the state-space to $S_m^{(2)} = \{2, 3, ..., m, 1\}$ and do this successively finishing up with $S_m^{(m)} = \{m, 1, 2, ... m-1\}$. This can be effected by permuting the elements of the transition matrix. For example, for $S_m^{(2)}$ we can do this by moving the elements of first column of $P^{(m)}$ to after the $m$-th column, followed by moving the first row to the last row, to obtain a new transition matrix $P^{(m)(2)}$. One of the easier ways to program this in *MATLAB* is to note that $P^{(m)(2)} (mod(row + m - 2, m) + 1, mod(col + m - 2, m) + 1) = P^{(m)(1)}(row, col)$.

Step 2: For $k = 2, 3, 4,..., m-1, m$.
(*i*) Repeat Step 1(*i*) with $P^{(m)} = P^{(m)(k)}$.
(*ii*) Repeat Step 1(*ii*) with $\boldsymbol{\mu}^{(k)(m)} = \boldsymbol{\mu}^{(m)} = (1,,1,...,1)$.
(*iii*) Repeat Step 1(*iii*) to calculate the $m \times 1$ column vector $\overline{\boldsymbol{m}}_m^{(k)(m)T} = \left(m_{kk}, m_{k+1,k},...,m_{k,m}, m_{k,1},...,m_{k-1,k}\right)$.

Step 3: Combine the results of the Steps 1(*iii*) and 2(*iii*) to find $M$ as follows.
Let $\overline{M} = (\boldsymbol{m}_m^{(1)(m)}, \overline{\boldsymbol{m}}_m^{(2)(m)},...,\overline{\boldsymbol{m}}_m^{(m)(m)})$ and reorder the elements of $\overline{M}$ to obtain $M = (\boldsymbol{m}_m^{(1)(m)}, \boldsymbol{m}_m^{(2)(m)},...,\boldsymbol{m}_m^{(m)(m)})$. This can be carried out in *MATLAB* by noting that for each row and column entry, $\overline{M}(mod(row + col - 2, m) + 1, col) = M(row, col)$.

A key observation is that the *EGTH* algorithm retains calculation accuracy since no subtractions are involved. Further, the stationary probabilities do not need to be computed in advance and can be found directly as inverses of the mean recurrence times.

## 5. Computation of *MFPTs* using modifications of the Heyman *FUND* algorithm

In carrying out Step 1 of the *EGTH* algorithm observe that the elements for $p_{ij}^{(n-1)}$ in the block upper left hand $(n-1) \times (n-1)$ corner of the transition matrix are based only on the elements $p_{ij}^{(n)}, p_{in}^{(n)}, p_{nj}^{(n)}$, $1 \leq i \leq n-1$, $1 \leq j \leq n-1$. This means that we can in effect overwrite the elements of the transition matrix that are not required in the future. At the conclusion of the reduction process we are left with a matrix of elements $\overline{P} = \left[\overline{p}_{ij}\right]$,

where 
$$\overline{p}_{ij} = \begin{cases} p_{ij}^{(j)} = u_{ij}, & 1 \leq i < j \leq m, \\ p_{ii}^{(i)} = d_{ii}, & 1 \leq i = j \leq m, \\ p_{ij}^{(i)} = l_{ij}, & 1 \leq j < i \leq m; \end{cases}$$
(5.1)

so that $\overline{P} = \overline{U} + \overline{D} + \overline{L}$

where $\overline{U} = [u_{ij}]$ is a strictly upper triangular matrix with zeros on and below the diagonal, $\overline{L} = [l_{ij}]$ is strictly lower triangular matrix with zeros on and above the diagonal and $\overline{D} = \text{diag}(d_{11}, \ldots, d_{mm})$ is a diagonal matrix.

From (4.2),

$$p_{ij}^{(n-2)} = p_{ij}^{(n-1)} + \frac{p_{i,n-1}^{(n-1)} p_{n-1,j}^{(n-1)}}{S(n-1)} = p_{ij}^{(n)} + \frac{p_{in}^{(n)} p_{nj}^{(n)}}{S(n)} + \frac{p_{i,n-1}^{(n-1)} p_{n-1,j}^{(n-1)}}{S(n-1)}, \quad 1 \le i \le n-2, \, 1 \le j \le n-2.$$

It is easy to establish, by considering $t = n-3, \ldots, n-k$, that

$$p_{ij}^{(t)} = p_{ij}^{(n)} + \sum_{k=t+1}^{n} \frac{p_{ik}^{(k)} p_{kj}^{(k)}}{S(k)}, \quad 1 \le i \le t \le n-1, \, 1 \le j \le t \le n-1.$$

Since $p_{ij}^{(m)} = p_{ij}$ using the notation of (5.1), that if $q_{kk} \equiv 1/S(k)$, for $k = 2, \ldots, m$, then

$$p_{ij}^{(t)} = p_{ij} + \sum_{k=t+1}^{m} u_{ik} q_{kk} l_{kj} \text{ for } 1 \le i \le t \le m-1, \, 1 \le j \le t \le m-1, \text{ with } p_{ij}^{(t)} = p_{ij} \text{ for } i = m \text{ or }$$

$j = m$. Thus $\overline{P} = P + \overline{U}\overline{Q}\overline{L}$ where $\overline{Q} \equiv \text{diag}(q_{11}, q_{22}, \ldots, q_{mm})$. Note that at this stage $q_{11}$ can be arbitrarily defined. The first column and last row of $\overline{U}$ are empty and the first row and column of $\overline{L}$ are also empty. Further, since $S(k) = \sum_{j=1}^{k-1} p_{kj}^{(k)} = 1 - p_{kk}^{(k)}$,

$\overline{D} = \text{diag}(p_{11}^{(1)}, p_{22}^{(2)}, \ldots, p_{mm}^{(m)}) = \text{diag}(1, 1-S(2), \ldots, 1-S(m))$ implying $I - \overline{D} = \text{diag}(0, S(2), \ldots, S(m))$.

Let $\overline{S} = \text{diag}(1, S(2), \ldots, S(m)) = E_{11} + I - \overline{D}$ so that $\overline{D} = E_{11} + I - \overline{S}$. We define $q_{11} = 1$ so that $\overline{S}^{-1} = \overline{Q}$. From these results we establish the following theorem.

**Theorem 5.1:** For an irreducible transition matrix $P$, the Markovian kernel $I - P$ can be factored into a $UL$ form where $L$ is a lower triangular matrix and $U$ is an upper triangular matrix, i.e. $I - P = UL$.

In particular, if $\overline{P} = \overline{U} + \overline{D} + \overline{L}$ is the matrix of overwritten elements of $P$ from the *GTH* algorithm, $U = \overline{U}\overline{S}^{-1} - I$ and $L = \overline{L} - (I - \overline{D})$ where $\overline{S} = E_{11} + (I - \overline{D})$.

**Proof**: From the results above

$$I - P = I - \overline{P} + \overline{U}\overline{Q}\overline{L} = I - \overline{U} - \overline{D} - \overline{L} + \overline{U}\overline{S}^{-1}\overline{L} = I - \overline{D} - \overline{U} + (\overline{U}\overline{S}^{-1} - I)\overline{L}$$

i.e. $I - P = \overline{S} - E_{11} - \overline{U} + (\overline{U}\overline{S}^{-1} - I)\overline{L}$ since $\overline{S} - E_{11} = I - \overline{D}$.

Now $(\overline{U}\overline{S}^{-1} - I)(E_{11} - \overline{S}) = \overline{U}\overline{S}^{-1} E_{11} - \overline{U} - E_{11} + \overline{S} = \overline{S} - E_{11} - \overline{U}$

since $\overline{U}\overline{S}^{-1} E_{11} = \overline{U} \text{diag}(q_{11}, q_{22}, \ldots, q_{mm}) E_{11} = \overline{U} E_{11} = 0$ (since $u_{11} = 0$)

Thus $I - P = (\overline{U}\overline{S}^{-1} - I)(\overline{L} - \overline{S} + E_{11}) = UL$, where

$U \equiv \overline{U}\overline{S}^{-1} - I$ is upper triangular and $L \equiv \overline{L} - \overline{S} + E_{11} = \overline{L} - (I - \overline{D})$ is lower triangular.

□

Grassman [4] first explored an *UL* factorisation of $I - P$ based upon the *GTH* algorithm. A version of this *UL* factorisation was used by Heyman [7] to produce his *FUND* algorithm to compute $Z$, the fundamental matrix of irreducible Markov chains. The proof given above is modified, due to some arbitrariness in the choice of the $\overline{Q}$ matrix, through a possible choice of $q_{11}$. Our choice for $\overline{Q}$ and hence for $\overline{S}$ leads to $U$ having all the elements of its diagonal as $-1$ and the other elements strictly upper triangular. This leads to $U$ having determinant $(-1)^m$ and consequently implying the non-singularity of $U$. $L$ has all the elements of its first row as 0.

Heyman [7] uses the $UL$ factorisation to find an expression for Z. We incorporate the results of his Theorem 1 within our Theorem 5.2 below but embed and extend his results in the g-inverse setting with a formal non-constructive proof.

**Theorem 5.2**: Let $P$ be the transition matrix of an irreducible finite $MC$, $\pi^T$ its stationary probability vector and $\Pi = e\pi^T$.
 *(a)* If $X$ is any solution of
$$(I - P)X = I - \Pi, \qquad (5.2)$$
then $X$ is a one-condition g-inverse of $I - P$ and satisfies the property that
$$Xe = xe, \text{ where } x \text{ is a constant}. \qquad (5.3)$$
*(b)* If $X$ is a solution of (5.2) then $A^{\#}$, the group inverse of $I - P$, is given by
$$A^{\#} = (I - \Pi)X, \qquad (5.4)$$
*(c)* If $X$ is a solution of (5.2) then $Z$, the fundamental inverse of $I - P$, is given by
$$Z = \Pi + (I - \Pi)X. \qquad (5.5)$$
**Proof:**
*(a)* Observe that from (5.2) and (1.2), $(I - P)X(I - P) = (I - e\pi^T)(I - P) = I - P$, implying that $X$ is a one-condition g-inverse of $I - P$. Further, from (5.2) and (1.2),
$(I - P)Xe = e - e\pi^T e = 0$, implying that $Xe$ is a right eigenvector of $I - P$ and hence must be a multiple of $e$ leading to (5.3).
*(b)* From Theorem 6.3 of [11] or Corollary 4.6 of [13], if $G$ is any g-inverse of $I - P$, when $P$ is irreducible, then $(I - \Pi)G(I - \Pi) = A^{\#}$. Taking $G = X$ observe that
$(I - \Pi)X(I - \Pi) = (I - \Pi)X - (I - \Pi)Xe\pi^T = (I - \Pi)X - xe\pi^T + xe\pi^T e\pi^T = (I - \Pi)X$ leading to (5.4).
*(c)* $Z$, the fundamental matrix of $I - P$, is given by $Z = [I - P + \Pi]^{-1} = A^{\#} + \Pi = \Pi + (I - \Pi)X$, leading to (5.5).
□

With $I - P = UL$, eqn. (5.2) can be solved in steps. Let
$$LX = Y, \qquad (5.6)$$
implying, from eqn. (5.2), that
$$UY = I - \Pi. \qquad (5.7)$$
We first solve, from eqn. (5.7), $Y$, uniquely, by backward substitution. In *MATLAB* we use the procedure $Y = U\setminus(I - \Pi)$. Note that for all $j = 1, \ldots, m$, $y_j^{(r)T}e = \sum_{i=1}^{m} y_{ij} = 0$.
Further, from eqn. (5.6), since $e_1^T L = \mathbf{0}^T$, we have that $e_1^T Y = \mathbf{0}^T$, and we conclude that the first row of $Y$, $y_1^{(r)T} = (y_{11}, \ldots, y_{1m})$, consists of zero elements.
Since we may take any one-condition g-inverse of $I - P$, we may take the first row of $X$ as the zero vector. Thus we may partition $L$, $X$ and $Y$ in block form as

$$LX = \begin{bmatrix} 0 & \mathbf{0}^T \\ \mathbf{l}_1^{(c)} & L_1 \end{bmatrix} \begin{bmatrix} 0 & \mathbf{0}^T \\ \mathbf{x}_1^{(c)} & X_1 \end{bmatrix} = \begin{bmatrix} 0 & \mathbf{0}^T \\ \mathbf{y}_1^{(c)} & Y_1 \end{bmatrix} = Y,$$

implying that $L_1 \mathbf{x}_1^{(c)} = \mathbf{y}_1^{(c)}$ and $L_1 X_1 = Y_1$, or equivalently $L_1(\mathbf{x}_1^{(c)}, X_1) = (\mathbf{y}_1^{(c)}, Y_1)$.
Thus if $\hat{X} = (\mathbf{x}_1^{(c)}, X_1)$, $\hat{Y} = (\mathbf{y}_1^{(c)}, Y_1)$, then we need to solve $L_1 \hat{X} = \hat{Y}$.

**Proc 10: (Heyman *FUND* Algorithm for *M* using *Z*).**
1. Start with $P$ and use the *GTH* algorithm to compute $\pi^T$.
2. Use the decomposition of Theorem 5.1 finding $\bar{P}$ and hence $U$ and $L$.
3. Solve $UY = I - \Pi$, where $\Pi = e\pi^T$, by back substitution.
4. Solve $L_1 \hat{X} = \hat{Y}$, by forward substitution.
5. Let $X = \begin{bmatrix} \mathbf{0}^T \\ \hat{X} \end{bmatrix}$.
6. Compute $Z = \Pi + (I - \Pi)X$.
7. Compute $M = [I - Z + EZ_d]D$ where $D = (\Pi_d)^{-1}$.

Heyman's *FUND* algorithm for finding the *MFPT*'s can be modified by noting, using (5.4), that instead of computing $Z$ one can compute the group inverse, reducing the number of calculations required in Step 6 of Proc 10 as follows.

**Proc 11: (Modified Heyman *FUND* Algorithm for *M* using $A^{\#}$).**
1. Carry out steps 1 to 5 of Proc 10.
2. Compute $A^{\#} = (I - \Pi)X$.
3. Compute $M = [I - A^{\#} + EA_d^{\#}]D$ where $D = (\Pi_d)^{-1}$.

One doesn't need to compute either $Z$ or $A^{\#}$ since $X$ is a one-condition g-inverse with the property that $Xe = \mathbf{0}$, ($x = 0$ in Theorem 5.2(*a*), since $X$ is chosen to have the first row the zero vector) and thus from Theorem 2.2(*c*) the following simpler Proc 12 is justified. Note that Heyman also observed this computational benefit for finding $M$ in the final section of his paper, [7].

**Proc 12: (Modified Heyman *FUND* Algorithm for *M* using *X*).**
1. Carry out steps 1 to 5 of Proc 10.
2. Compute $M = [I - X + EX_d]D$ where $D = (\Pi_d)^{-1}$.

**6. Test problems**

We use the following test problems that were introduced by Harrod & Plemmons [5]. as poorly conditioned examples for computing the stationary distribution of the underlying irreducible *MC* as well as having been used as examples for testing various different algorithms for computing *M*, the matrix of *MFPTs*, ([8], [9]). While the dimensions of the state-space are relatively small, the test problems lead to some computational difficulties.

**TP1:** (As modified by Heyman and Reeves ([9]). The original version of *TP1*, given in [5] related to a 10-state MC however it was shown, by Heyman [6], that four of the states were in fact transient and the irreducible sub chain was identified as

$$\begin{bmatrix} .1 & .6 & 0 & .3 & 0 & 0 \\ .5 & .5 & 0 & 0 & 0 & 0 \\ .5 & .2 & 0 & 0 & .3 & 0 \\ 0 & .7 & 0 & .2 & 0 & .1 \\ .1 & 0 & .8 & 0 & 0 & .1 \\ .4 & 0 & .4 & 0 & 0 & .2 \end{bmatrix}$$

**TP2:** A typo for the original problem for element (1,5) was identified and corrected in [9]. The test problem is also known as the 8 X 8 Courtois matrix and was also considered in a paper by Benzi [1].

$$\begin{bmatrix} .85 & 0 & .149 & .0009 & 0 & .00005 & 0 & .00005 \\ .1 & .65 & .249 & 0 & .0009 & .00005 & 0 & .00005 \\ .1 & .8 & .09996 & .0003 & 0 & 0 & .0001 & 0 \\ 0 & .0004 & 0 & .7 & .2995 & 0 & .0001 & 0 \\ .0005 & 0 & .0004 & .399 & .6 & .0001 & 0 & 0 \\ 0 & .00005 & 0 & 0 & .00005 & .6 & .2499 & .15 \\ .00003 & 0 & .00003 & .00004 & 0 & .1 & .8 & .0999 \\ 0 & .00005 & 0 & 0 & .00005 & .1999 & .25 & .55 \end{bmatrix}.$$

**TP3:**

$$\begin{bmatrix} 0.999999 & 1.0\,E-07 & 2.0\,E-07 & 3.0\,E-07 & 4.0\,E-07 \\ 0.4 & 0.3 & 0 & 0 & 0.3 \\ 5.0\,E-07 & 0 & 0.999999 & 0 & 5.0\,E-07 \\ 5.0\,E-07 & 0 & 0 & 0.999999 & 5.0\,E-07 \\ 2.0\,E-07 & 3.0\,E-07 & 1.0\,E-07 & 4.0\,E-07 & 0.999999 \end{bmatrix}.$$

**TP4 and variants:**
**TP41** $:\varepsilon =1.0E-01;$ **TP42** $:\varepsilon =1.0E-03;$ **TP43** $:\varepsilon =1.0E-05;$ **TP44** $:\varepsilon =1.0E-07$

$$\begin{bmatrix} .1-\varepsilon & .3 & .1 & .2 & .3 & \varepsilon & 0 & 0 & 0 & 0 \\ .2 & .1 & .1 & .2 & .4 & 0 & 0 & 0 & 0 & 0 \\ .1 & .2 & .2 & .4 & .1 & 0 & 0 & 0 & 0 & 0 \\ .4 & .2 & .1 & .2 & .1 & 0 & 0 & 0 & 0 & 0 \\ .6 & .3 & 0 & 0 & .1 & 0 & 0 & 0 & 0 & 0 \\ \varepsilon & 0 & 0 & 0 & 0 & .1-\varepsilon & .2 & .2 & .4 & .1 \\ 0 & 0 & 0 & 0 & 0 & .2 & .2 & .1 & .3 & .2 \\ 0 & 0 & 0 & 0 & 0 & .1 & .5 & 0 & .2 & .2 \\ 0 & 0 & 0 & 0 & 0 & .5 & .2 & .1 & 0 & .2 \\ 0 & 0 & 0 & 0 & 0 & .1 & .2 & .2 & .3 & .2 \end{bmatrix}$$

**TP5** and **TP6:** The previous examples have small state-spaces. We are also interested in how the algorithms perform for larger state-spaces. As we could not find any specific test cases in the literature we include two cases where we generate the transition matrices for a 100-state (*TP*5) and a 500-state (*TP*6) Markov chain. To have some sparsity we used *MATLAB* to generate $A = rand(m) > 0.6$, an $m \times m$ matrix of zero's and one's with the long run proportion of zero elements being 0.6. We removed the diagonal elements to construct $B = A - diag(diag(A))$, followed by computing $TP = inv(diag(sum(B')))*B$. We established, in each case, that the minimum entry of the square of *TP* is positive so that *TP* is the transition matrix of an irreducible *MC*. The transition matrices were stored in separate .mat files to ensure that the same *TP* was used for each procedure. The mat files for *TP*5 and *TP*6 are available from the author.

## 7. Computational error comparisons

For numerical computations and comparisons, we coded each algorithm using *MATLAB* (64-bit version R2015b on a *MacBook Air* computer) and used the test problems of Section 6. *MATLAB* was run in both single and double precision to enable us to compute and compare the *MFPT* matrices $M(S) = [m_{ij}(S)]$ and $M(D) = [m_{ij}(D)]$.

Before discussing the merits of the different procedures we make some general comments. Proc 1 appears to suffer from the requirement to compute the specific matrix inverse that, according to [8], can lead to significant inaccuracy, (with three types of errors – computing the stationary probabilities, constructing the matrix and the computation of the inverse.) Using the *GTH* algorithm minimises the first type of error. The matrix leads to a dense matrix even, if one starts with a sparse matrix, and the inverse is potentially fraught with computational difficulties for large state-space cases. Proc 2 may also have matrix computation difficulties for large state-spaces but it does eliminate the need to find the $\pi_j$ probabilities and does preserve sparsity to a limited extent. As mentioned earlier, we do not expect these two procedures to perform well in large state-space cases. Procs 3 to 8 are all based on perturbation procedures that, because of their coding simplicity are typically computationally efficient. However it is not entirely clear how well the accuracy at each step is maintained. Proc 3 is based on an update of g-inverses. Proc 4 involves the update of the group inverse (that may introduce an additional unnecessary calculations). Proc 5 extends the number of computations required requiring two sequences of matrix updates, while Procs 6, 7 and 8, all have a simple internal structure for updating rank one modifications with different initial and final conditions. One would expect these latter three procedures to have similar performance characteristics. Proc 9 is based on state reduction and enlargement and has the desirable feature that no subtraction operation need be performed. Since the procedure cycles amongst all starting states no decision need be taken in respect to the renumbering states to determine an "ideal" starting state. The way we have constructed the algorithm may however take longer computation times as we are using only the *MFPT*'s calculation in a single column at any iteration. Proc 10 is based on the idea of Heyman using the fundamental matrix as does Proc 11 using the group inverse but these introduce additional computations which are simplified in Proc 12. Based on the above observations, our expectation is that Procs 9 and 12 should be among the better performing procedures for accuracy.

We are interested in determining whether we can determine if one or some of the procedures lead to more accurate results for the elements of *M* than others.

In the first instance we ran all the procedures using *MATLAB* and reviewed the output for each *MFPT* matrix *M*. All procedures ran without problem except for *TP*44 under single precision. A detailed examination of the output for this test problem, under single precision, show that warnings are expressed in *MATLAB* for Procedures 1 and 2 indicating that "the matrix is close to singular or badly scaled. Results may be inaccurate". For Proc 1, eight *MFPT* elements of *M* are zero. Close observation of the Proc 2 output indicates that some of the small *MFPT*s are very inaccurate. For Proc 3, fifty elements of *M* are negative, which we did not expect. For Proc 4 the message "Warning: Matrix is singular to working precision" is given with the net effect that no terms can be calculated. All of the four Procs 5, 6, 7 and 8 each yield the elements $m_{21}$, $m_{31}$, $m_{41}$ and $m_{51}$ being negative. In Procs 10 and 11, $m_{76} < 0$, $m_{86} < 0$, $m_{10,6} < 0$ and $m_{9,6} = 0$, while in Proc 12 $m_{76} = 0$ and $m_{86} = 0$. The only procedure to come through unscathed is Proc 9. There were no problems with the larger state-space cases *TP*5 and *TP*6 except for larger errors. These observations suggest that single precision calculations should, where possible, be avoided for all the procedures.

Consider the "errors" $\varepsilon_{ij} = m_{ij} - \sum_{k \neq j} p_{ik} m_{kj} - 1$, where the $m_{ij}$ have been calculated by the relevant algorithmic procedures. By virtue of Eqn. (1.3), ideally these errors $\varepsilon_{ij}$ would all be zero for a perfect computation. In this paper we explore three sensible measures:

1)  *PZE(D)*, be "*the percentage of error terms $\varepsilon_{ij}$ that are zero*", under double precision.

2)  $ORE(D) = \sum_{i=1}^{m} \sum_{j=1}^{m} |\varepsilon_{ij}|$, the "*overall residual error*", computed in double precision.

3)  $ANED = \frac{1}{m^2} \sum_{i=1}^{m} \sum_{j=1}^{m} \left\{ -\log_{10} \left| \frac{m_{ij}(D) - m_{ij}(S)}{m_{ij}(D)} \right| \right\}$, t*he average number of extra digits*

achieved by the double precision calculation $m_{ij}(D)$ over the single precision calculation $m_{ij}(S)$.

There are a number of other general measures that we can use to explore the computation errors, both under single precision and double precision. We omit these as they typically lead to similar conclusions to those deduced in this paper. They are included in the latest ArXiv.com version of this paper, (arXiv:1701.0778).

**Percentage of Zero errors:**
If a high proportion of the $m^2$ error terms are zero for a procedure we would expect such a procedure to be have a high level of accuracy. For all the perturbation Procedures 3 to 8 for *TP*2, *TP*3, *TP*42, *TP*43 and *TP*44, *PZE(S) (*under single precision*)* is 0% and less than 0.75% for *TP*5 and 0.12% for *TP*6. This immediately suggests that single precision calculations are not all sensible for the perturbation procedures.

Using double precision for *PZE(D)*, (see Table 1 in the Appendix), Proc 9 has the highest proportion of zero errors for each small state-space test problem except for *TP*1, when it is outperformed by Proc 12. For the larger state-space cases, Proc 2 (*TP*5) and Proc 7 (*TP*6) gave the largest values for *PZE(D)*, while for these cases, Procedures 3, 8 and 9 all had a smaller number of zero errors than the other procedures that were typically in the range 18.01 to 18.90 for *TP*5 and 16.94 to 17.51 for *TP*6.

In general, the worse performing procedures are Proc 3 (*TP*1, *TP*41, *TP*42, *TP*44, *TP*5), Proc 4 (*TP*2) and Proc 8 (*TP*3, *TP*43, *TP*6).

**Overall Residual errors:**
The *ORE* is one of the better indicators of accuracy. Minimum values of *ORE(D)* are achieved by Proc 9 for all test problems in double precision, except for *TP*1, *TP*5 and *TP* 6 where the minimum values are achieved by Proc 11 (*TP*1) and Proc 1 (*TP*5, *TP*6). Note that Proc 1 is the well-known standard procedure of Kemeny and Snell – still a reliable computational tool that involves a matrix inversion. Its accuracy is surprising since it requires the initial computation of the stationary distribution (admittedly, using the very accurate *GTH* algorithm.) It is also interesting to note that Proc 2, a relatively simple matrix form for computing the *MFPTs* that also involves a matrix inverse, performs well in just about all situations. These results are due to a very efficient matrix inversion procedure incorporated in *MATLAB* and is contrary to the comments expressed by Heyman and O'Leary [8] where they state that "Deriving means … of first passage times from … the fundamental matrix Z … leads to significant inaccuracy on the more difficult problems."

The worse performing procedures in respect to *ORE(D)* are Proc 3 (*TP*1, *TP*41, *TP*42, *TP*44, *TP*5), Proc 4 (*TP*2, *TP*3), Proc 8 (*TP*43, *TP*6).

The perturbation procedures, and in particular Procs 3, 4, and 5 do not perform well in general. Procs 6, 7 and 8 are all variants of a perturbation procedure and do not star in general, except for Proc 6 in *TP*2, (4[th] best), *TP*3 (2[nd] best), *TP*42 (3[rd] best), *TP*43 (4[th] best), *TP*44 (3[rd] best), *TP*5 (3[rd] best) and *TP*6 (2[nd] best). The variants of Heyman's FUND algorithm, Procs 10, 11 and 12 all perform well for all test problems except for *TP*3.

It is interesting to note that the favored algorithm, Proc 9, does not perform as well expected in the larger state-space cases, (*TP*5 and *TP*6).

**The average number of extra digits for double precision over single precsion:**
The ANED statistic was introduced by Heyman and Reeves [9] and Heyman and O'Leary [8], in comparing *MFPT* calculations, where one regards the double precision result as the "true" result and the single precision result as the "computed" result, taking the *number of (extra) accurate digits)* can be defined as the overall average of

$$\log_{10} \left| \frac{result_{true} - result_{computed}}{result_{true}} \right|.$$

In both [8] and [9] the results for *ANED* were displayed in figures and no actual numerical results were tabulated.

Note that *TP*1 has some unique features in respect to computing *ANED*. We actually deduce exact results for three *MFPT's* since it can be shown that $m_{21} = 2$, $m_{43} = 160.5$ and $m_{53} = 26.3$. In computing the *MFPT* matrices, under double precision, all twelve procedures obtain these exact three results. Under single precision, only Proc 1 and Proc 12 yield all three exact results, while Proc 2 gives the exact results for $m_{21}$ and $m_{43}$, and Procedures 9, 10 and 11 yield the exact result for only $m_{21}$. Thus when calculating the average number of accurate digits we must omit the results when the *MFPTs* under single and double precision are the same, as the logarithm of zero is negative infinity. In Table 3 of the Appendix for *ANED* for *TP*1, we indicate with *** when the average is taken over the 33 finite terms, ** with an average over 34 terms and * over 35 terms.

Proc 9 gives the largest *ANED* statistic for all the test problems except for *TP*5 and *TP*6 when Proc 8, the best, followed by Procs 5, 6 and 7 all perform better than Proc 9. The worse performing procedure for all *TP*s is Proc 3 (even with a negative value for *TP*44.) Proc 4 also performs uniformly poorly.

Heyman and O'Leary [8] used two algorithms for computing *M*, without previously computing *Z* or $A^{\#}$ – an algorithm *MH*, which is based on the two stage *UL* factorisation, using $U^{-1}$ and *L* factors as in Heyman's method (similar to our Proc 12), and an algorithm *M*, using the *UL* factors and normalisation. They obtained values, for the number of extra accurate digits, between 6 and 7 for all *TPs* for their algorithm *M* but displayed widely varying values for different *TPs* for the algorithm *MH*.

Heyman and Reeves [9] presented four algorithms - *LINPACK*, *SR*, *KSGTH*, and *KSGAUSS* for computing *M* with the same test problems used in this paper. They explored different software packages deciding that *LINPACK* "worked the best" although the solution computed by *LINPACK* did not run for *TP*44 as the matrix inverse could not be computed. The *KSGTH* is the same as our Proc 1 using the *GTH* algorithm to compute the steady state probabilities while *KSGAUSS* used Gaussian elimination to solve the stationary equations. Their favoured algorithm, is the *SR*, "State reduction", procedure of Kohlas, on which Proc 9 is based. The Hunter *EGTH* Procedure consistently produces results for *ANED* in the range 7.30 to 7.43, for the small state space cases, similar to that achieved by Heyman & Reeves [9] (as extrapolated from their graphical output). Of interest is the *ANED* values of Proc 8 for the large state space cases, viz, 8.30 for *TP*5 and 9.09 for *TP*6.

There are a range of other error comparisons that we can make but they generally end up pointing to Proc 9, Hunter's *EGTH* Procedure, as giving the most accurate results for the small state cases. Our general recommendation is to use this procedure. While there was slight drop in precision for the larger state-space cases, Proc 8 which performed the best in these two cases did not figure prominently in the smaller state space cases. It is clear that single precision is not generally recommended as a suitable computation procedure as it is fraught with inaccuracies.

When paper [18] was written it was hoped that the perturbation procedures were going to generally yield very accurate results, comparable with other procedures. Apart from isolated situations, Procedures 3 – 8 do not perform as well as we had hoped. The *EGTH* procedure, involving no subtractions, completely overshadows the perturbation procedures in the small state cases. Proc 12, while not in the same class as Proc 9, reliably produces the second most accurate results in the small state-space cases.

## 8. Computational time comparisons

The computation time of the procedures is a secondary consideration in this paper as often programming codes and computing platforms can results in considerable variation. However we used the *MATLAB* 'tic-toc' procedure, with the same version (*R*2015*b*) and computer (*MacBook Air*) as used for the error computations. Each algorithm was run 10 times and the elapsed computing times were averaged. The results are displayed in Table 4 in the Appendix. Under double precision Proc 2 gave the fastest times for all *TP*s except for *TP*2, *TP*3 and *TP*41 when Proc 5 was the fastest. The slowest procedures were Proc 6 for all *TP* except Proc 9 for the large state cases, *TP*5 and *TP*6. We knew that while the coding of Proc 9 was simple, we were in effect discarding the computation of many of the *MFPTs* produced at each iteration. The algorithm as designed actually contains all the *MFPTs* in columns 1 and 2 of the *MFPT* matrix at each run so we are actually discarding elements that have been computed. Further, the rearrangement of the transition matrix at each iteration adds significantly to the computational times. Some skilful redesigning based on the theoretical results presented in Hunter [17] could potentially lead to a considerable reduction in computational time. This, of course, would detract from the overall structural simplicity of the algorithm that really comes into its own if one is simply interested in computing just the first column of *M* to find the $m_{i1}$, $i = 1, …, n$ elements.

## 9. Conclusions

As discussed earlier, when this article was initially conceived it was the accuracy of the computations that was of paramount importance, with the author under the impression that modern computing techniques would take care of any time considerations. At the suggestion of a referee two larger state-space cases and computation times for each of the procedures were included in this final version. Those two suggestions highlighted that small-state *MC* examples and large state-case *MC* examples can sometimes behave differently both in accuracy and in relative computation times.

In summary, for small state-space cases, the top four performing procedures (in recommended order of accuracy) are Procs 9, 12, 2 and 6. Proc 2 has the advantage of being computationally the fastest with Proc 12 also relatively fast to implement.

For large state-space cases, the top four performing procedures (in recommended order of accuracy) are Procs 1, 6, 2 and 12. Three of these four procedures (Proc 2, 1, 12) are also faster to implement than any others.

In both small-state and large-state situations Procs 3 and 8 are not at all useful. We also do not recommend using Procs 4 and 5. The only perturbation procedure that we recommend for accuracy is Proc 6 (due in part to its simple initial conditions) but this is, however, the slowest procedure for the small-state cases.

The disappointment is that the best small-state case performer Proc 9 does not appear to perform well for large-state cases.

**Acknowledgements**
The author would like to express his thanks to the referee who suggested including some examples of larger, possibly sparse, state-space Markov chains in the paper, as well as exploring the times to carry out the computations. The author sought to find some examples of the transition matrices of some larger state-space Markov chains in the open literature but was unable to track down anything suitable. He wishes to acknowledge the suggestions made by the Associate Editor, Professor Steve Kirkland that led to the construction of *TP*5 and *TP*6.

**Appendix 1: Selected Error calculations for all Procedures and all Test Problems**

Table 1: Percentage of Zero Error Terms, under double precision

| PZE(D) | TP1 | TP2 | TP3 | TP41 | TP42 | TP43 | TP44 | TP5 | TP6 |
|---|---|---|---|---|---|---|---|---|---|
| Proc 1 | 47.22 | 9.38 | 20.00 | 24.00 | 19.00 | 17.00 | 25.00 | 18.86 | 17.45 |
| Proc 2 | 33.33 | 21.88 | 24.00 | 21.00 | 23.00 | 22.00 | 32.00 | 18.90 | 17.38 |
| Proc 3 | 5.56 | 7.81 | 12.00 | 7.00 | 7.00 | 12.00 | 9.00 | 7.10 | 3.45 |
| Proc 4 | 25.00 | 4.69 | 8.00 | 14.00 | 9.00 | 22.00 | 19.00 | 18.79 | 17.34 |
| Proc 5 | 36.11 | 7.81 | 12.00 | 10.00 | 17.00 | 13.00 | 18.00 | 18.20 | 16.94 |
| Proc 6 | 52.78 | 17.19 | 16.00 | 31.00 | 23.00 | 23.00 | 25.00 | 18.34 | 17.39 |
| Proc 7 | 38.89 | 14.06 | 8.00 | 28.00 | 22.00 | 17.00 | 19.00 | 18.64 | 17.51 |
| Proc 8 | 27.78 | 7.81 | 4.00 | 17.00 | 11.00 | 11.00 | 12.00 | 7.74 | 2.95 |
| Proc 9 | 47.22 | 23.44 | 28.00 | 34.00 | 27.00 | 31.00 | 34.00 | 11.85 | 7.07 |
| Proc 10 | 41.67 | 20.31 | 16.00 | 28.00 | 23.00 | 22.00 | 23.00 | 18.29 | 17.29 |
| Proc 11 | 44.44 | 20.31 | 16.00 | 32.00 | 24.00 | 20.00 | 23.00 | 18.01 | 17.40 |
| Proc 12 | 55.56 | 23.44 | 16.00 | 33.00 | 25.00 | 20.00 | 31.00 | 18.52 | 17.43 |

Table 2: Overall Residual Error, under double precision

| ORE(D) | TP1 | TP2 | TP3 | TP41 | TP42 | TP43 | TP44 | TP5 | TP6 |
|---|---|---|---|---|---|---|---|---|---|
| Proc 1 | 3.149E-13 | 4.294E-11 | 8.240E-09 | 4.347E-13 | 4.159E-11 | 3.435E-09 | 3.290E-07 | 1.460E-10 | 2.017E-08 |
| Proc 2 | 2.981E-13 | 3.865E-11 | 7.215E-09 | 4.408E-13 | 2.837E-11 | 4.093E-09 | 1.969E-07 | 1.491E-10 | 2.028E-08 |
| Proc 3 | 1.108E-12 | 1.091E-10 | 1.660E-08 | 1.626E-12 | 1.007E-10 | 8.704E-09 | 9.925E-07 | 4.325E-10 | 1.165E-07 |
| Proc 4 | 4.083E-13 | 3.527E-10 | 1.635E-04 | 6.856E-13 | 6.806E-11 | 4.751E-09 | 5.257E-07 | 1.539E-10 | 2.047E-08 |
| Proc 5 | 5.854E-13 | 8.020E-11 | 1.478E-08 | 7.641E-13 | 6.991E-11 | 6.965E-09 | 7.254E-07 | 1.783E-10 | 2.236E-08 |
| Proc 6 | 3.132E-13 | 3.964E-11 | 6.906E-09 | 3.804E-13 | 2.939E-11 | 3.231E-09 | 2.283E-07 | 1.501E-10 | 2.022E-08 |
| Proc 7 | 3.632E-13 | 4.894E-11 | 1.138E-08 | 3.685E-13 | 4.454E-11 | 4.660E-09 | 4.237E-07 | 1.507E-10 | 2.038E-08 |
| Proc 8 | 4.582E-13 | 1.044E-10 | 1.571E-08 | 8.028E-13 | 8.569E-11 | 1.177E-08 | 6.241E-07 | 3.897E-10 | 1.176E-07 |
| Proc 9 | 2.955E-13 | 2.848E-11 | 5.275E-09 | 2.883E-13 | 1.956E-11 | 1.577E-09 | 1.417E-07 | 2.580E-10 | 5.901E-08 |
| Proc 10 | 2.449E-13 | 4.127E-11 | 5.211E-05 | 3.447E-13 | 3.794E-11 | 3.129E-09 | 3.452E-07 | 1.540E-10 | 2.037E-08 |
| Proc 11 | 2.169E-13 | 4.217E-11 | 5.211E-05 | 3.223E-13 | 3.749E-11 | 3.303E-09 | 3.452E-07 | 1.545E-10 | 2.038E-08 |
| Proc 12 | 2.849E-13 | 3.643E-11 | 5.211E-05 | 2.970E-13 | 3.071E-11 | 3.048E-09 | 2.566E-07 | 1.519E-10 | 2.027E-08 |

Table 3: Average number of extra digits, of double precision over single precision

| ANED | TP1 | TP2 | TP3 | TP41 | TP42 | TP43 | TP44 | TP 5 | TP6 |
|---|---|---|---|---|---|---|---|---|---|
| Proc 1 | ***6.752 | 4.969 | 3.225 | 6.509 | 5.042 | 3.227 | 0.782 | 6.884 | 6.714 |
| Proc 2 | **7.118 | 4.513 | 2.459 | 6.743 | 5.100 | 3.563 | 2.311 | 6.904 | 6.866 |
| Proc 3 | 6.277 | 3.960 | 1.742 | 6.037 | 3.988 | 2.169 | -0.198 | 6.625 | 6.025 |
| Proc 4 | 6.381 | 4.515 | 1.895 | 6.344 | 4.098 | 2.200 | NaN | 7.045 | 6.853 |
| Proc 5 | 6.821 | 4.346 | 2.101 | 6.565 | 4.525 | 2.510 | 0.107 | 8.766 | 8.957 |
| Proc 6 | 6.778 | 4.353 | 2.101 | 6.768 | 4.803 | 2.787 | 0.502 | 7.532 | 7.321 |
| Proc 7 | 6.801 | 4.149 | 2.101 | 6.493 | 4.474 | 2.453 | 0.145 | 7.353 | 7.224 |
| Proc 8 | 6.908 | 4.353 | 2.101 | 6.740 | 4.803 | 2.787 | 0.502 | 8.827 | 9.088 |
| Proc 9 | *7.350 | 7.293 | 7.353 | 7.368 | 7.416 | 7.430 | 7.332 | 7.073 | 6.810 |
| Proc 10 | *6.6272 | 4.509 | 2.968 | 6.875 | 5.172 | 3.312 | 1.479 | 6.825 | 6.489 |
| Proc 11 | *6.595 | 4.509 | 2.968 | 6.934 | 5.172 | 3.312 | 1.479 | 6.827 | 6.485 |
| Proc 12 | ***6.766 | 4.544 | 2.968 | 7.043 | 5.583 | 4.159 | 2.740 | 6.913 | 6.538 |

Table 4: Average computation times (seconds) for M in double precision

| Mdouble | TP1 | TP2 | TP3 | TP41 | TP42 | TP43 | TP44 | TP5 | TP6 |
|---|---|---|---|---|---|---|---|---|---|
| Proc 1 | 0.0329 | 0.0426 | 0.0415 | 0.0402 | 0.0384 | 0.0388 | 0.0364 | 0.4950 | 3.3186 |
| Proc 2 | 0.0149 | 0.0209 | 0.0202 | 0.0212 | 0.0183 | 0.0198 | 0.0194 | 0.0262 | 0.0843 |
| Proc 3 | 0.0188 | 0.0216 | 0.0197 | 0.0212 | 0.0211 | 0.0206 | 0.0206 | 4.4360 | 35.0946 |
| Proc 4 | 0.0289 | 0.0359 | 0.0257 | 0.0270 | 0.0272 | 0.0271 | 0.0273 | 5.3691 | 42.1983 |
| Proc 5 | 0.0194 | 0.0201 | 0.0197 | 0.0210 | 0.0259 | 0.0199 | 0.0206 | 6.2228 | 46.5435 |
| Proc 6 | 0.1174 | 0.1167 | 0.1093 | 0.1143 | 0.1206 | 0.1149 | 0.1084 | 1.8839 | 16.7556 |
| Proc 7 | 0.0429 | 0.0438 | 0.0419 | 0.0434 | 0.0429 | 0.0444 | 0.0459 | 1.6789 | 15.8758 |
| Proc 8 | 0.0301 | 0.0302 | 0.0289 | 0.0307 | 0.0319 | 0.0312 | 0.0300 | 1.7948 | 15.8933 |
| Proc 9 | 0.0388 | 0.0404 | 0.0417 | 0.0450 | 0.0441 | 0.0419 | 0.0425 | 168.3714 | 1156.2318 |
| Proc 10 | 0.0748 | 0.0762 | 0.0723 | 0.0792 | 0.0804 | 0.0770 | 0.0770 | 0.5351 | 3.4063 |
| Proc 11 | 0.0792 | 0.0753 | 0.0877 | 0.0751 | 0.0805 | 0.0771 | 0.0769 | 0.5345 | 3.4054 |
| Proc 12 | 0.0740 | 0.0715 | 0.0683 | 0.0708 | 0.0772 | 0.0941 | 0.0764 | 0.5320 | 3.3882 |

**Appendix 2: Additional Error calculations**

In the version that appears in "Linear Algebra and its Applications" (available at https://doi.org/10.1016/j.laa.2018.03.010 ) three measures of accuracy of the computations were considered, *viz*. *PZE(D), ORE(D)* and *ANED*.

From the errors associated with each *(i,j)* element in computing the MFPT $m_{ij}$ from state *i* to state *j*, $\varepsilon_{ij} = m_{ij} - \sum_{k \neq j} p_{ik} m_{kj} - 1$, where the $m_{ij}$ have been calculated by the relevant algorithmic procedures, we can compute a variety of additional measures to explore the computation errors, both under single precision and double precision. In particular, the *minimum absolute residual errors*, $MINARE(.) = \min_{1 \leq i \leq m, 1 \leq j \leq m} |\varepsilon_{ij}|$, and the *maximum absolute residual error*, $MAXARE(.) = \max_{1 \leq i \leq m, 1 \leq j \leq m} |\varepsilon_{ij}|$,

as well as *ORE(S)*, the single precision version of the overall residual error.

We can also compare the single precision and double precision calculations in terms of the *minimum absolute error*, $MINAE(S, D) = \min_{1 \leq i \leq m, 1 \leq j \leq m} |m_{ij}(S) - m_{ij}(D)|$,

*the maximum absolute error*, $MAXAE(S, D) = \max_{1 \leq i \leq m, 1 \leq j \leq m} |m_{ij}(S) - m_{ij}(D)|$,

*and the relative error between the double and single precision computations as*
$$REL(S, D) = \sum_{i=1}^{m} \sum_{j=1}^{m} |m_{ij}(S) - m_{ij}(D)|.$$

We can also consider *PCZE(S)*, the single precision version of "*the percentage of error terms $\varepsilon_{ij}$ that are zero*".

The *MATLAB* code for computation of these measures are given Appendix 3.
Tables of these measures are given in Appendix 4 and some charts are provided in Appendix 5.

**Appendix 3: MATLAB Code for all the procedures and error calculations**

The input variables are *m*, *TP.mat* (which gives the transition matrix *P*). (Alternatively one can list the entries of P in standard MATLAB format)

The outputs are the matrices MS and MD, the single precision and double precision versions of the mean first passage time matrices.

Use PART 1 of the procedure, in single precision, to find MS and calculate the single precision errors for the procedure:

```
deltaMS= MS-PS*MS-ES+PS*DS;
MINARES=min(min(abs(deltaMS)))
MAXARES= max(max(abs(deltaMS)))
ORES=sum(sum(abs(deltaMS)))
idxS=deltaMS==0;
outS=sum(idxS(:));
PCZES=outS*100/(m*m)
```

Then use PART 2 of the procedure, in double precision, to find MD and calculate the double precision errors for the procedure:

```
deltaMD= MD-PD* MD-ED+PD*DD;
MINARED=min(min(abs(deltaMD)))
MAXARED= max(max(abs(deltaMD)))
ORED=sum(sum(abs(deltaMD)))
idxD=deltaMD==0;
outD=sum(idxD(:));
PCZED=outD*100/(m*m)
```

Then calculate the single and double comparison errors:

```
absMSD=abs(MS-MD);
MINAESD=min(min(absMSD))
MAXAESD=max(max(absMSD))
RELSD=sum(sum(absMSD))
```

Finally, compute the ANED statistic
```
A= MD;
B= MS;
D=abs(A-B);
for i=1:m;
j=1:m;
K(i,j)=D(i,j)./A(i,j);
end;
L=-log10(K);
G=sum(sum(L));
ANED = G/(m*m)
```

**Procedure 1**
**PART 1**
```
echo off
clear all
format long
m=;
load('TP.mat');
PS=single(P);
eS=[eye(m)];
e1S=ones(m,1);
ES=ones(m,m);
IS=eye(m);
PPS=PS;
AASS=zeros(m,m);
for n=m:-1:2
   SS(1,n)=sum(PPS(n,1:n-1));
     for i=1:n-1
      for j=1:n-1
       AASS(i,j)=PPS(i,n)*PPS(n,j)/SS(1,n);
       PPS(i,j)=PPS(i,j)+AASS(i,j);
       end
     end
end
rS=zeros(1,m);
rS(1,1)=1;
for n=2:m
    for i=1:n-1
     rS(1,n)=rS(1,n)+rS(1,i)*PPS(i,n)/SS(1,n);
     end
end
TOTS=sum(rS);
pit_1S=rS/TOTS;
PiS=e1S*pit_1S;
ZS=inv(IS-PS+PiS);
DS=inv(diag(diag(PiS)));
MS=(IS-ZS+ES*diag(diag(ZS)))*DS;
```

**PART 2**
```
PD=double(P);
eD=[eye(m)];
e1D=ones(m,1);
ED=ones(m,m);
ID=eye(m);
PPD=PD;
AASD=zeros(m,m);
for n=m:-1:2
   SD(1,n)=sum(PPD(n,1:n-1));
     for i=1:n-1
      for j=1:n-1
       AASD(i,j)=PPD(i,n)*PPD(n,j)/SD(1,n);
       PPD(i,j)=PPD(i,j)+AASD(i,j);
       end
     end
end
rD=zeros(1,m);
rD(1,1)=1;
for n=2:m
    for i=1:n-1
     rD(1,n)=rD(1,n)+rD(1,i)*PPD(i,n)/SD(1,n);
     end
end
TOTD=sum(rD);
pit_1D=rD/TOTD;
PiD=e1D*pit_1D;
ZD=inv(ID-PD+PiD);
DD=inv(diag(diag(PiD)));
MD=(ID-ZD+ED*diag(diag(ZD)))*DD;
```

**Procedure 2**
**PART 1**
```
echo off
clear all
format long
m=;
```

```matlab
load('TP.mat')
PS=single(P);
eS=[eye(m)];
e1S=ones(m,1);
ES=ones(m,m);
IS=eye(m);
gS=inv(IS-PS+e1S*eS(:,1)');
pit_S=gS(1,:);
MSS=ones(m,m);
for i=1:m
   for j=1:m
      MSS(i,j)=(gS(j,j)-gS(i,j))/gS(1,j);
   end
   MSS(i,i)=1/gS(1,i);
end
DS=inv(diag(diag(e1S*pit_S)));
MS=MSS;
```

**PART 2**
```matlab
PD=double(P);
eD=[eye(m)];
e1D=ones(m,1);
ED=ones(m,m);
ID=eye(m);
gD=inv(ID-PD+e1D*eD(:,1)');
pit_D=gD(1,:);
MSD=ones(m,m);
for i=1:m
   for j=1:m
      MSD(i,j)=(gD(j,j)-gD(i,j))/gD(1,j);
   end
   MSD(i,i)=1/gD(1,i);
end
DD=inv(diag(diag(e1D*pit_D)));
MD=MSD;
```

**Procedure 3**
**PART 1**
```matlab
echo off
clear all
format long
m=;
load('TP.mat')
PS=single(P);
eS=[ones(m,1) eye(m)];
e0S=ones(m,1);
ES=ones(m,m);
US=[ones(m,1)/m zeros(m,m)];
GS=zeros(m,m,m+1);
GS(:,:,1)=eye(m);
I=eye(m);
PiS=zeros(m,m);
gS=zeros(1,m);
HS=zeros(m,m,m);
DSS=zeros(m,m,m);
MSS=zeros(m,m,m);
for i=1:m
   US(:,i+1)=US(:,i)+PS(i,:)'-e0S/m;
   gS(1,i)=US(:,i)'*GS(:,:,i)*eS(:,1+i);
   GS(:,:,i+1)=GS(:,:,i)+GS(:,:,i)*(eS(:,i)-eS(:,i+1))*(US(:,i)'*GS(:,:,i)/gS(1,i));
   PiS(:,i)=(US(:,i+1)'*GS(:,:,i+1)/(US(:,i+1)'*GS(:,:,i+1)*e0S))';
   HS(:,:,i)=GS(:,:,i+1)*(I-e0S*PiS(:,i)');
   DSS(:,:,i)=inv(diag(diag(e0S*PiS(:,i)')));
   MSS(:,:,i)=(I-HS(:,:,i)+ES*(diag(diag(HS(:,:,i)))))*DSS(:,:,i);
end
DS=DSS(:,:,m);
MS=MSS (:,:,m);
```

**PART 2**
```matlab
PD=double(P);
eD=[ones(m,1) eye(m)];
e0D=ones(m,1);
ED=ones(m,m);
UD=[ones(m,1)/m zeros(m,m)];
```

```
GD=zeros(m,m,m+1);
GD(:,:,1)=eye(m);
I=eye(m);
PiD=zeros(m,m);
gD=zeros(1,m);
HD=zeros(m,m,m);
DDD=zeros(m,m,m);
MDD=zeros(m,m,m);
for i=1:m
   UD(:,i+1)=UD(:,i)+PD(i,:)'-e0D/m;
   gD(1,i)=UD(:,i)'*GD(:,:,i)*eD(:,1+i);
   GD(:,:,i+1)=GD(:,:,i)+GD(:,:,i)*(eD(:,i)-eD(:,i+1))*(UD(:,i)'*GD(:,:,i)/gD(1,i));
   PiD(:,i)=(UD(:,i+1)'*GD(:,:,i+1)/(UD(:,i+1)'*GD(:,:,i+1)*e0D))';
    HD(:,:,i)=GD(:,:,i+1)*(I-e0D*PiD(:,i)');
    DDD(:,:,i)=inv(diag(diag(e0D*PiD(:,i)')));
    MDD(:,:,i)=(I-HD(:,:,i)+ED*(diag(diag(HD(:,:,i)))))*DDD(:,:,i);
end
PitD=PiD(:,m)';
DD=DDD(:,:,m);
MD=MDD(:,:,m);
```

**Procedure 4**
**PART 1**
```
echo off
clear all
format long
m=;
load('TP.mat')
PS=single(P);
eS=[ones(m,1) eye(m)];
e0S=ones(m,1);
ES=ones(m,m);
IS=eye(m);
P0S=ones(m,m)/m;
kS=zeros(1,m);
RS=zeros(m,m,m+1);
RS(:,:,1)=IS-P0S;
PiS=zeros(m,m);
bS=zeros(m,m);
DSS=zeros(m,m,m);
AsharpS=zeros(m,m,m+1);
AsharpS(:,:,1)=IS-P0S;
MSS=zeros(m,m,m);
for i=1:m
   bS(:,i)=PS(i,:)'-e0S/m;
   kS(1,i)=1-bS(:,i)'*RS(:,:,i)*eS(:,i+1);
   RS(:,:,i+1)=RS(:,:,i)+1/kS(i)*RS(:,:,i)*eS(:,i+1)* bS(:,i)'*RS(:,:,i);
   PiS(:,i)=(eS(:,2)-(eS(:,2)'*(IS-PS)*RS(:,:,i+1))')';
   AsharpS(:,:,i+1)=(IS-e0S*PiS(:,i)')*RS(:,:,i+1);
   DSS(:,:,i)=inv(diag(diag(e0S*PiS(:,i)')));
   MSS(:,:,i)=(IS-AsharpS(:,:,i+1)+ES*(diag(diag(AsharpS(:,:,i+1)))))*DSS(:,:,i);
end
DS=DSS(:,:,m);
MS=MSS(:,:,m);
```

**PART 2**
```
PD=double(P);
eD=[ones(m,1) eye(m)];
e0D=ones(m,1);
ED=ones(m,m);
ID=eye(m);
P0D=ones(m,m)/m;
kD=zeros(1,m);
RD=zeros(m,m,m+1);
RD(:,:,1)=ID-P0D;
PiD=zeros(m,m);
bD=zeros(m,m);
DD=zeros(m,m,m);
AsharpD=zeros(m,m,m+1);
AsharpD(:,:,1)=ID-P0D;
MDD=zeros(m,m,m);
for i=1:m
   bD(:,i)=PD(i,:)'-e0D/m;
   kD(1,i)=1-bD(:,i)'*RD(:,:,i)*eD(:,i+1);
```

```
    RD(:,:,i+1)=RD(:,:,i)+1/kD(i)*RD(:,:,i)*eD(:,i+1)* bD(:,i)'*RD(:,:,i);
    PiD(:,i)=(eD(:,2)-(eD(:,2)'*(ID-PD)*RD(:,:,i+1))')';
   AsharpD(:,:,i+1)=(ID-e0D*PiD(:,i)')*RD(:,:,i+1);
    DD(:,:,i)=inv(diag(diag(e0D*PiD(:,i)')));
    MDD(:,:,i)=(ID-AsharpD(:,:,i+1)+ES*(diag(diag(AsharpD(:,:,i+1)))))*DD(:,:,i);
end
DD= DD(:,:,m);
MD=MDD (:,:,m);
```

## Procedure 5
### PART 1
```
echo off
clear all
format long
m=;
load('TP.mat');
PS=single(P);
eS=[ones(m,1) eye(m)];
e0S=ones(m,1);
ES=ones(m,m);
IS=eye(m);
P0S=ones(m,m)/m;
AsharpS=zeros(m,m,m+1);
AsharpS(:,:,1)=IS-P0S;
fS=zeros(1,m);
SS=zeros(m,m,m);
PiS=zeros(m,m,m+1);
PiS(:,:,1)=P0S;
bS=zeros(m,m);
DSS=zeros(m,m,m);
MSS=zeros(m,m,m);
for i=1:m
    bS(:,i)=PS(i,:)'-e0S/m;
    fS(1,i)=1-bS(:,i)'*AsharpS(:,:,i)*eS(:,i+1);
    SS(:,:,i)=IS+1/fS(i)*eS(:,i+1)* bS(:,i)'*AsharpS(:,:,i);
    PiS(:,:,i+1)=PiS(:,:,i)*SS(:,:,i);
    AsharpS(:,:,i+1)=(IS-PiS(:,:,i+1))*AsharpS(:,:,i)*SS(:,:,i);
    DSS(:,:,i)=inv(diag(diag(PiS(:,:,i+1))));
    MSS(:,:,i)=(IS-AsharpS(:,:,i+1)+ES*(diag(diag(AsharpS(:,:,i+1)))))*DSS(:,:,i);
end
DS=DSS (:,:,m);
MS=MSS (:,:,m);
```

### PART 2
```
PD=double(P);
P0D=ones(m,m)/m;
e0D=ones(m,1);
eD=[ones(m,1) eye(m)];
ED=ones(m,m);
ID=eye(m);
AsharpD=zeros(m,m,m+1);
AsharpD(:,:,1)=ID-P0D;
fD=zeros(1,m);
SD=zeros(m,m,m);
PiD=zeros(m,m,m+1);
PiD(:,:,1)=P0D;
bD=zeros(m,m);
DDD=zeros(m,m,m);
MDD=zeros(m,m,m);
for i=1:m
    bD(:,i)=PD(i,:)'-e0D/m;
    fD(1,i)=1-bD(:,i)'*AsharpD(:,:,i)*eD(:,i+1);
    SD(:,:,i)=ID+1/fD(i)*eD(:,i+1)* bD(:,i)'*AsharpD(:,:,i);
    PiD(:,:,i+1)=PiD(:,:,i)*SD(:,:,i);
    AsharpD(:,:,i+1)=(ID-PiD(:,:,i+1))*AsharpD(:,:,i)*SD(:,:,i);
    DDD(:,:,i)=inv(diag(diag(PiD(:,:,i+1))));
    MDD(:,:,i)=(ID-AsharpD(:,:,i+1)+ED*(diag(diag(AsharpD(:,:,i+1)))))*DDD(:,:,i);
end
DD= DDD(:,:,m);
MD=MDD (:,:,m);
```

## Procedure 6
### PART 1
```
echo off
```

```
clear all
format long
m=;
load('TP.mat')
PS=single(P);
eS=[eye(m)];
e0S=ones(m,1);
ES=ones(m,m);
US=[ones(m,1)/m zeros(m,m)];
GS=zeros(m,m,m+1);
GS(:,:,1)=eye(m);
IS=eye(m);
btS=ones(m,m);
for i=1:m
    btS(i,:)=PS(i,:)-1/m*e0S';
end
KS=ones(m,m,m+1);
kS=ones(m,1);
CS=ones(m,m,m);
KS(:,:,1)=IS;
for i=1:m
    kS(i,1)=1-btS(i,:)*KS(:,:,i)*eS(:,i);
    CS(:,:,i)=1/kS(i,1)*eS(:,i)*btS(i,:)*KS(:,:,i);
    KS(:,:,i+1)=KS(:,:,i)*(IS+CS(:,:,i));
end
KKS=KS(:,:,m+1);
pitS=1/m*e0S'*KKS;
DS=inv(diag(diag(e0S*pitS)));
MS=(IS-KKS+ES*diag(diag(KKS)))*DS;
```

**PART 2**
```
PD=double(P);
eD=[eye(m)];
e0D=ones(m,1);
ED=ones(m,m);
UD=[ones(m,1)/m zeros(m,m)];
GD=zeros(m,m,m+1);
GD(:,:,1)=eye(m);
ID=eye(m);
btD=ones(m,m);
for i=1:m
    btD(i,:)=PD(i,:)-1/m*e0D';
end
KD=ones(m,m,m+1);
kD=ones(m,1);
CD=ones(m,m,m);
KD(:,:,1)=ID;
for i=1:m
    kD(i,1)=1-btD(i,:)*KD(:,:,i)*eD(:,i);
    CD(:,:,i)=1/kD(i,1)*eD(:,i)*btD(i,:)*KD(:,:,i);
    KD(:,:,i+1)=KD(:,:,i)*(ID+CD(:,:,i));
end
KKD=KD(:,:,m+1);
pitD=1/m*e0D'*KKD;
DD=inv(diag(diag(e0D*pitD)));
MD=(ID-KKD+ED*diag(diag(KKD)))*DD;
```

**Procedure 7**
**PART 1**
```
echo off
clear all
format long
m=;
load('TP.mat')
PS=single(P);
eS=[eye(m)];
e0S=ones(m,1);
ES=ones(m,m);
US=[ones(m,1)/m zeros(m,m)];
GS=zeros(m,m,m+1);
GS(:,:,1)=eye(m);
IS=eye(m);
btS=ones(m,m);
for i=1:m
```

```
   btS(i,:)=PS(i,:)-1/m*e0S';
end
KS=ones(m,m,m+1);
kS=ones(m,1);
CS=ones(m,m,m);
KS(:,:,1)=IS+e0S*(e0S'/m-eS(:,1)');
for i=1:m
   kS(i,1)=1-btS(i,:)*KS(:,:,i)*eS(:,i);
   CS(:,:,i)=1/kS(i,1)*eS(:,i)*btS(i,:)*KS(:,:,i);
   KS(:,:,i+1)=KS(:,:,i)*(IS+CS(:,:,i));
end
KKS=KS(:,:,m+1);
pitS=eS(:,1)'*KKS;
DS=inv(diag(diag(e0S*pitS)));
MS=(IS-KKS+ES*diag(diag(KKS)))*DS;
```

**PART 2**
```
PD=double(P);
eD=[eye(m)];
e0D=ones(m,1);
ED=ones(m,m);
UD=[ones(m,1)/m zeros(m,m)];
GD=zeros(m,m,m+1);
GD(:,:,1)=eye(m);
ID=eye(m);
btD=ones(m,m);

for i=1:m
   btD(i,:)=PD(i,:)-1/m*e0D';
end
KD=ones(m,m,m+1);
kD=ones(m,1);
CD=ones(m,m,m);
KD(:,:,1)=ID+ e0D*(e0D'/m-eD(:,1)');
for i=1:m
   kD(i,1)=1-btD(i,:)*KD(:,:,i)*eD(:,i);
   CD(:,:,i)=1/kD(i,1)*eD(:,i)*btD(i,:)*KD(:,:,i);
   KD(:,:,i+1)=KD(:,:,i)*(ID+CD(:,:,i));
end
KKD=KD(:,:,m+1);
pitD=eD(:,1)'*KKD;
DD=inv(diag(diag(e0D*pitD)));
MD=(ID-KKD+ED*diag(diag(KKD)))*DD;
```

**Procedure 8**
**PART 1**
```
echo off
clear all
format long
m=;
load('TP.mat')
PS=single(P);
eS=[eye(m)];
e0S=ones(m,1);
ES=ones(m,m);
US=[ones(m,1)/m zeros(m,m)];
GS=zeros(m,m,m+1);
GS(:,:,1)=eye(m);
IS=eye(m);
btS=ones(m,m);
for i=1:m
   btS(i,:)=PS(i,:)-1/m*e0S';
end
KS=ones(m,m,m+1);
kS=ones(m,1);
CS=ones(m,m,m);
KS(:,:,1)=IS-(m-1)/m^2*e0S*e0S';
for i=1:m
   kS(i,1)=1-btS(i,:)*KS(:,:,i)*eS(:,i);
   CS(:,:,i)=1/kS(i,1)*eS(:,i)*btS(i,:)*KS(:,:,i);
   KS(:,:,i+1)=KS(:,:,i)*(IS+CS(:,:,i));
end
KKS=KS(:,:,m+1);
pitS=e0S'*KKS;
```

```matlab
DS=inv(diag(e0S*pitS));
MS=(IS-KKS+ES*diag(diag(KKS)))*DS;
```

**PART 2**
```matlab
PD=double(P);
eD=[eye(m)];
e0D=ones(m,1);
ED=ones(m,m);
UD=[ones(m,1)/m zeros(m,m)];
GD=zeros(m,m,m+1);
GD(:,:,1)=eye(m);
ID=eye(m);
btD=ones(m,m);
for i=1:m
    btD(i,:)=PD(i,:)-1/m*e0D';
end
KD=ones(m,m,m+1);
kD=ones(m,1);
CD=ones(m,m,m);
KD(:,:,1)=ID-(m-1)/m^2*e0D*e0D';
for i=1:m
    kD(i,1)=1-btD(i,:)*KD(:,:,i)*eD(:,i);
    CD(:,:,i)=1/kD(i,1)*eD(:,i)*btD(i,:)*KD(:,:,i);
    KD(:,:,i+1)=KD(:,:,i)*(ID+CD(:,:,i));
end
KKD=KD(:,:,m+1);
pitD=e0D'*KKD;
DD=inv(diag(e0D*pitD));
MD=(ID-KKD+ED*diag(diag(KKD)))*DD;
```

**Procedure 9**
**PART 1**
```matlab
echo off
clear all
format long
m=;
load('TP.mat')
TMs= single(P);
Ps=TMs;
PPs=TMs;
es=ones(m,1);
ets= ones(1,m);
Ss=ones(1,m);
ES=ones(m,m);
mus=zeros(m,m);
mus(:,m)=1;
Ms=zeros(m,m);
for k=1:m
    for n=m:-1:2
    Ss(1,n)=sum(PPs(n,1:n-1));
    for i=1:n-1
      for j=1:n-1
         PPs(i,j)=PPs(i,j)+PPs(i,n)*PPs(n,j)/Ss(1,n);
      end
      mus(i,n-1)=mus(i,n)+mus(n,n)*PPs(i,n)/Ss(1,n);
    end
    end
    Ms(1,k)=(PPs(2,1)*mus(1,2)+PPs(1,2)*mus(2,2))/PPs(2,1);
    for n=2:m
      mms=0;
      for i=2:n-1
      mms=mms+PPs(n,i)*Ms(i,k);
      end
      Ms(n,k)=(mms+mus(n,n))/Ss(1,n);
    end
    for col=1:m
      for row= 1:m
      P_news1(mod(row+m-2,m)+1,col)=Ps(row,col);
      end
    end
    for col=1:m
      for row= 1:m
      P_news2(row,mod(col+m-2,m)+1)=P_news1(row,col);
      end
```

```
        end
     Ps=P_news2;
     PPs=Ps;
   end
   for col=1:m
     for row=1:m
        M_GTHs(mod(row+col-2,m)+1,col)=Ms(row,col);
     end
   end
M_GTHs;
PS=TMs;
DS=diag(diag(M_GTHs));
MS= M_GTHs;
```

**PART 2**
```
TMd= double(P);
Pd=TMd;
PPd=TMd;
ed=ones(m,1);
etd= ones(1,m);
Sd=ones(1,m);
ED=ones(m,m);
mud=zeros(m,m);
mud(:,m)=1;
Md=zeros(m,m);
for k=1:m
     for n=m:-1:2
     Sd(1,n)=sum(PPd(n,1:n-1));
     for i=1:n-1
       for j=1:n-1
          PPd(i,j)=PPd(i,j)+PPd(i,n)*PPd(n,j)/Sd(1,n);
       end
       mud(i,n-1)=mud(i,n)+mud(n,n)*PPd(i,n)/Sd(1,n);
     end;
     end;
     Md(1,k)=(PPd(2,1)*mud(1,2)+PPd(1,2)*mud(2,2))/PPd(2,1);
     for n=2:m
        mmd=0;
        for i=2:n-1
        mmd=mmd+PPd(n,i)*Md(i,k);
        end
        Md(n,k)=(mmd+mud(n,n))/Sd(1,n);
     end;
     for col=1:m
       for row= 1:m
       P_newd1(mod(row+m-2,m)+1,col)=Pd(row,col);
        end;
     end
     for col=1:m
       for row= 1:m
       P_newd2(row,mod(col+m-2,m)+1)=P_newd1(row,col);
        end
     end
     Pd=P_newd2;
     PPd=Pd;
  end
  for col=1:m
     for row=1:m
        M_GTHd(mod(row+col-2,m)+1,col)=Md(row,col);
     end
  end
M_GTHd;
MD= M_GTHd
DD=diag(diag(MD));
PD=TMd
```

**Procedure 10**
**PART 1**
```
echo off
clear all
format long
m=;
load('TP.mat')
PS=single(P);
```

```
PPS=single(P);
eS=ones(m,1);
ES=ones(m,m);
IS=eye(m);
SS=ones(1,m);
AAS=zeros(m,m);
for n=m:-1:2
   SS(1,n)=sum(PPS(n,1:n-1));
   for i=1:n-1
      for j=1:n-1
      AAS(i,j)=PPS(i,n)*PPS(n,j)/SS(1,n);
      PPS(i,j)=PPS(i,j)+AAS(i,j);
      end
   end
end
rS=zeros(1,m);
rS(1,1)=1;
for n=2:m
    for i=1:n-1
     rS(1,n)=rS(1,n)+rS(1,i)*PPS(i,n)/SS(1,n);
     end
end
TOTS=sum(rS);
PiGTHS=rS'/TOTS;
PbarS=PPS(1:m,1:m);
LbarS=tril(PbarS,-1);
UbarS=triu(PbarS,1);
DbarS=diag(diag(PbarS));
E11S= transpose(IS(1,1:m))* IS(1,1:m);
SbarS=IS+E11S-DbarS;
US=UbarS*inv(SbarS)-IS;
LS=LbarS-IS+DbarS;
ULfactS=US*LS;
eS=ones(m,1);
WS= eS*transpose(PiGTHS);
BS=IS-WS;
YS=US\BS;
L1S=LS(2:m,2:m);
Y1S=YS(2:m,1:m);
X1S=L1S\Y1S;
XS=[zeros(1,m);X1S];
ZS=WS+BS*XS;
DS=inv(diag(diag(WS)));
MS=(IS-ZS+ES*diag(diag(ZS)))*DS;
```

**PART 2**
```
PD=double(P);
PPD=double(P);
eD=ones(m,1);
ED=ones(m,m);
ID=eye(m);
SD=ones(1,m);
AAD=zeros(m,m);
for n=m:-1:2
   SD(1,n)=sum(PPD(n,1:n-1));
   for i=1:n-1
      for j=1:n-1
      AAD(i,j)=PPD(i,n)*PPD(n,j)/SD(1,n);
      PPD(i,j)=PPD(i,j)+AAD(i,j);
      end
   end
end
rD=zeros(1,m);
rD(1,1)=1;
for n=2:m
    for i=1:n-1
     rD(1,n)=rD(1,n)+rD(1,i)*PPD(i,n)/SD(1,n);
     end
end
TOTD=sum(rD);
PiGTHD=rD'/TOTD;
PbarD=PPD(1:m,1:m);
LbarD=tril(PbarD,-1);
UbarD=triu(PbarD,1);
```

```
DbarD=diag(diag(PbarD));
E11D= transpose(ID(1,1:m))* ID(1,1:m);
SbarD=ID+E11D-DbarD;
UD=UbarD*inv(SbarD)-ID;
LD=LbarD-ID+DbarD;
ULfactD=UD*LD;
eD=ones(m,1);
WD= eD*transpose(PiGTHD);
BD=ID-WD;
YD=UD\BD;
L1D=LD(2:m,2:m);
Y1D=YD(2:m,1:m);
X1D=L1D\Y1D;
XD=[zeros(1,m);X1D];
ZD=WD+BD*XD;
DD=inv(diag(diag(WD)));
MD=(ID-ZD+ED*diag(diag(ZD)))*DD;
```

**Procedure 11**:
**PART 1**
```
echo off
clear all
format long
m=;
load('TP.mat')
PS=single(P);
PPS=single(P);
eS=ones(m,1);
ES=ones(m,m);
IS=eye(m);
SS=ones(1,m);
AAS=zeros(m,m);
for n=m:-1:2
   SS(1,n)=sum(PPS(n,1:n-1));
   for i=1:n-1
      for j=1:n-1
      AAS(i,j)=PPS(i,n)*PPS(n,j)/SS(1,n);
      PPS(i,j)=PPS(i,j)+AAS(i,j);
      end
   end
end
rS=zeros(1,m);
rS(1,1)=1;
for n=2:m
    for i=1:n-1
     rS(1,n)=rS(1,n)+rS(1,i)*PPS(i,n)/SS(1,n);
    end
end
TOTS=sum(rS);
PiGTHS=rS'/TOTS;
PbarS=PPS(1:m,1:m);
LbarS=tril(PbarS,-1);
UbarS=triu(PbarS,1);
DbarS=diag(diag(PbarS));
E11S= transpose(IS(1,1:m))* IS(1,1:m);
SbarS=IS+E11S-DbarS;
US=UbarS*inv(SbarS)-IS;
LS=LbarS-IS+DbarS;
ULfactS=US*LS;
eS=ones(m,1);
WS= eS*transpose(PiGTHS);
BS=IS-WS;
YS=US\BS;
L1S=LS(2:m,2:m);
Y1S=YS(2:m,1:m);
X1S=L1S\Y1S;
XS=[zeros(1,m);X1S];
AsharpS=BS*XS;
DS=inv(diag(diag(WS)));
MS=(IS-AsharpS+ES*diag(diag(AsharpS)))*DS;
```

**PART 2**
```
PD=double(P);
PPD=double(P);
```

```
eD=ones(m,1);
ED=ones(m,m);
ID=eye(m);
SD=ones(1,m);
AAD=zeros(m,m);
for n=m:-1:2
   SD(1,n)=sum(PPD(n,1:n-1));
   for i=1:n-1
      for j=1:n-1
      AAD(i,j)=PPD(i,n)*PPD(n,j)/SD(1,n);
      PPD(i,j)=PPD(i,j)+AAD(i,j);
      end
   end
end
rD=zeros(1,m);
rD(1,1)=1;
for n=2:m
    for i=1:n-1
      rD(1,n)=rD(1,n)+rD(1,i)*PPD(i,n)/SD(1,n);
    end
end
TOTD=sum(rD);
PiGTHD=rD'/TOTD;
PbarD=PPD(1:m,1:m);
LbarD=tril(PbarD,-1);
UbarD=triu(PbarD,1);
DbarD=diag(diag(PbarD));
E11D= transpose(ID(1,1:m))* ID(1,1:m);
SbarD=ID+E11D-DbarD;
UD=UbarD*inv(SbarD)-ID;
LD=LbarD-ID+DbarD;
ULfactD=UD*LD;
eD=ones(m,1);
WD= eD*transpose(PiGTHD);
BD=ID-WD;
YD=UD\BD;
L1D=LD(2:m,2:m);
Y1D=YD(2:m,1:m);
X1D=L1D\Y1D;
XD=[zeros(1,m);X1D];
AsharpD=BD*XD;
DD=inv(diag(diag(WD)));
MD=(ID-AsharpD+ED*diag(diag(AsharpD)))*DD;
```

**Procedure 12**
**PART 1**
```
echo off
clear all
format long
m=;
load('TP.mat')
PS=single(P);
PPS=single(P);
eS=ones(m,1);
ES=ones(m,m);
IS=eye(m);
SS=ones(1,m);
AAS=zeros(m,m);
for n=m:-1:2
   SS(1,n)=sum(PPS(n,1:n-1));
   for i=1:n-1
      for j=1:n-1
      AAS(i,j)=PPS(i,n)*PPS(n,j)/SS(1,n);
      PPS(i,j)=PPS(i,j)+AAS(i,j);
      end
   end
end
rS=zeros(1,m);
rS(1,1)=1;
for n=2:m
    for i=1:n-1
      rS(1,n)=rS(1,n)+rS(1,i)*PPS(i,n)/SS(1,n);
    end
end
```

```
TOTS=sum(rS);
PiGTHS=rS'/TOTS;
PbarS=PPS(1:m,1:m);
LbarS=tril(PbarS,-1);
UbarS=triu(PbarS,1);
DbarS=diag(diag(PbarS));
E11S= transpose(IS(1,1:m))* IS(1,1:m);
SbarS=IS+E11S-DbarS;
US=UbarS*inv(SbarS)-IS;
LS=LbarS-IS+DbarS;
ULfactS=US*LS;
eS=ones(m,1);
WS= eS*transpose(PiGTHS);
BS=IS-WS;
YS=US\BS;
L1S=LS(2:m,2:m);
Y1S=YS(2:m,1:m);
X1S=L1S\Y1S;
XS=[zeros(1,m);X1S];
DS=inv(diag(diag(WS)));
MS=(IS-XS+ES*diag(diag(XS)))*DS;
```

**PART 2**
```
PD=double(P);
PPD=double(P);
eD=ones(m,1);
ED=ones(m,m);
ID=eye(m);
SD=ones(1,m);
AAD=zeros(m,m);
for n=m:-1:2
   SD(1,n)=sum(PPD(n,1:n-1));
   for i=1:n-1
      for j=1:n-1
      AAD(i,j)=PPD(i,n)*PPD(n,j)/SD(1,n);
      PPD(i,j)=PPD(i,j)+AAD(i,j);
      end
   end
end
rD=zeros(1,m);
rD(1,1)=1;
for n=2:m
    for i=1:n-1
      rD(1,n)=rD(1,n)+rD(1,i)*PPD(i,n)/SD(1,n);
    end
end
TOTD=sum(rD);
PiGTHD=rD'/TOTD;
PbarD=PPD(1:m,1:m);
LbarD=tril(PbarD,-1);
UbarD=triu(PbarD,1);
DbarD=diag(diag(PbarD));
E11D= transpose(ID(1,1:m))* ID(1,1:m);
SbarD=ID+E11D-DbarD;
UD=UbarD*inv(SbarD)-ID;
LD=LbarD-ID+DbarD;
ULfactD=UD*LD;
eD=ones(m,1);
WD= eD*transpose(PiGTHD);
BD=ID-WD;
YD=UD\BD;
L1D=LD(2:m,2:m);
Y1D=YD(2:m,1:m);
X1D=L1D\Y1D;
XD=[zeros(1,m);X1D];
DD=inv(diag(diag(WD)));
MD=(ID-XD+ED*diag(diag(XD)))*DD;
```

# Appendix 4: Error calculations for all Procedures and all Test Problems

## Table 4.1: MINARE under Single Precision

| MINARE(S) | TP1 | TP2 | TP3 | TP41 | TP42 | TP43 | TP44 | TP5 | TP6 |
|---|---|---|---|---|---|---|---|---|---|
| Proc 1 | 0.0000E+00 | 0.0000E+00 | 0.0000E+00 | 0.0000E+00 | 0.0000E+00 | 0.0000E+00 | 1.0902E-02 | 0.0000E+00 | 0.0000E+00 |
| Proc 2 | 0.0000E+00 | 0.0000E+00 | 0.0000E+00 | 0.0000E+00 | 0.0000E+00 | 0.0000E+00 | 0.0000E+00 | 0.0000E+00 | 0.0000E+00 |
| Proc 3 | 0.0000E+00 | 2.8133E-05 | 1.5333E-02 | 1.1921E-07 | 3.3379E-06 | 1.9091E-04 | 1.9091E-04 | 0.0000E+00 | 0.0000E+00 |
| Proc 4 | 0.0000E+00 | 2.8014E-06 | 1.2846E-03 | 0.0000E+00 | 2.5034E-06 | 1.4853E-04 | NaN | 0.0000E+00 | 0.0000E+00 |
| Proc 5 | 0.0000E+00 | 9.9784E-06 | 1.2996E-03 | 0.0000E+00 | 1.1921E-06 | 1.8179E-04 | 4.9167E-02 | 0.0000E+00 | 0.0000E+00 |
| Proc 6 | 0.0000E+00 | 1.0908E-05 | 1.2996E-03 | 0.0000E+00 | 5.8413E-06 | 8.6933E-04 | 1.0224E-01 | 0.0000E+00 | 0.0000E+00 |
| Proc 7 | 0.0000E+00 | 5.7817E-06 | 1.2996E-03 | 0.0000E+00 | 1.5497E-06 | 3.1388E-04 | 4.2600E+00 | 0.0000E+00 | 0.0000E+00 |
| Proc 8 | 0.0000E+00 | 1.7941E-05 | 1.2996E-03 | 0.0000E+00 | 4.4107E-06 | 2.2805E-04 | 7.9652E-02 | 0.0000E+00 | 0.0000E+00 |
| Proc 9 | 0.0000E+00 | 0.0000E+00 | 0.0000E+00 | 0.0000E+00 | 0.0000E+00 | 0.0000E+00 | 0.0000E+00 | 0.0000E+00 | 0.0000E+00 |
| Proc 10 | 0.0000E+00 | 0.0000E+00 | 0.0000E+00 | 0.0000E+00 | 0.0000E+00 | 0.0000E+00 | 8.6164E-03 | 0.0000E+00 | 0.0000E+00 |
| Prov 11 | 0.0000E+00 | 0.0000E+00 | 0.0000E+00 | 0.0000E+00 | 0.0000E+00 | 0.0000E+00 | 8.6164E-03 | 0.0000E+00 | 0.0000E+00 |
| Proc 12 | 0.0000E+00 | 0.0000E+00 | 0.0000E+00 | 0.0000E+00 | 0.0000E+00 | 0.0000E+00 | 0.0000E+00 | 0.0000E+00 | 0.0000E+00 |

## Table 4.2: MINARE under Double Precision

| MINARE(D) | TP1 | TP2 | TP3 | TP41 | TP42 | TP43 | TP44 | TP5 | TP6 |
|---|---|---|---|---|---|---|---|---|---|
| Proc 1 | 0.0000E+00 | 0.0000E+00 | 0.0000E+00 | 0.0000E+00 | 0.0000E+00 | 0.0000E+00 | 0.0000E+00 | 0.0000E+00 | 0.0000E+00 |
| Proc 2 | 0.0000E+00 | 0.0000E+00 | 0.0000E+00 | 0.0000E+00 | 0.0000E+00 | 0.0000E+00 | 0.0000E+00 | 0.0000E+00 | 0.0000E+00 |
| Proc 3 | 0.0000E+00 | 0.0000E+00 | 0.0000E+00 | 0.0000E+00 | 0.0000E+00 | 0.0000E+00 | 0.0000E+00 | 0.0000E+00 | 0.0000E+00 |
| Proc 4 | 0.0000E+00 | 0.0000E+00 | 0.0000E+00 | 0.0000E+00 | 0.0000E+00 | 0.0000E+00 | 0.0000E+00 | 0.0000E+00 | 0.0000E+00 |
| Proc 5 | 0.0000E+00 | 0.0000E+00 | 0.0000E+00 | 0.0000E+00 | 0.0000E+00 | 0.0000E+00 | 0.0000E+00 | 0.0000E+00 | 0.0000E+00 |
| Proc 6 | 0.0000E+00 | 0.0000E+00 | 0.0000E+00 | 0.0000E+00 | 0.0000E+00 | 0.0000E+00 | 0.0000E+00 | 0.0000E+00 | 0.0000E+00 |
| Proc 7 | 0.0000E+00 | 0.0000E+00 | 0.0000E+00 | 0.0000E+00 | 0.0000E+00 | 0.0000E+00 | 0.0000E+00 | 0.0000E+00 | 0.0000E+00 |
| Proc 8 | 0.0000E+00 | 0.0000E+00 | 0.0000E+00 | 0.0000E+00 | 0.0000E+00 | 0.0000E+00 | 0.0000E+00 | 0.0000E+00 | 0.0000E+00 |
| Proc 9 | 0.0000E+00 | 0.0000E+00 | 0.0000E+00 | 0.0000E+00 | 0.0000E+00 | 0.0000E+00 | 0.0000E+00 | 0.0000E+00 | 0.0000E+00 |
| Proc 10 | 0.0000E+00 | 0.0000E+00 | 0.0000E+00 | 0.0000E+00 | 0.0000E+00 | 0.0000E+00 | 0.0000E+00 | 0.0000E+00 | 0.0000E+00 |
| Proc 11 | 0.0000E+00 | 0.0000E+00 | 0.0000E+00 | 0.0000E+00 | 0.0000E+00 | 0.0000E+00 | 0.0000E+00 | 0.0000E+00 | 0.0000E+00 |
| Proc 12 | 0.0000E+00 | 0.0000E+00 | 0.0000E+00 | 0.0000E+00 | 0.0000E+00 | 0.0000E+00 | 0.0000E+00 | 0.0000E+00 | 0.0000E+00 |

## Table 4.3: MAXARE under Single Precision

| MAXARE(S) | TP1 | TP2 | TP3 | TP41 | TP42 | TP43 | TP44 | TP5 | TP6 |
|---|---|---|---|---|---|---|---|---|---|
| Proc 1 | 1.5259E-05 | 2.3131E-03 | 1.0000E+00 | 7.6294E-06 | 7.3242E-04 | 9.3750E-02 | 1.9313E+01 | 4.5776E-05 | 3.0518E-04 |
| Proc 2 | 3.0518E-05 | 1.9531E-03 | 1.0000E+00 | 1.5259E-05 | 9.7656E-04 | 6.2500E-02 | 7.7565E+00 | 6.9618E-05 | 5.4932E-04 |
| Proc 3 | 5.8830E-05 | 1.9503E-03 | 7.4928E-01 | 9.0599E-06 | 8.2135E-04 | 5.8108E-02 | 5.8108E-02 | 4.1604E-05 | 2.5988E-04 |
| Proc 4 | 3.0518E-05 | 1.3188E-03 | 7.6065E-01 | 6.4373E-06 | 7.2205E-04 | 4.4751E-02 | NaN | 4.3154E-05 | 2.7853E-04 |
| Proc 5 | 5.7161E-05 | 1.1806E-03 | 1.4523E+00 | 5.8413E-06 | 3.8600E-04 | 5.0105E-02 | 2.7472E+01 | 3.8862E-05 | 2.9123E-04 |
| Proc 6 | 2.2531E-05 | 1.5189E-03 | 1.4884E+00 | 5.7220E-06 | 3.3975E-04 | 4.0965E-02 | 1.1226E+01 | 3.7432E-05 | 2.7871E-04 |
| Proc 7 | 8.1241E-05 | 1.9780E-03 | 1.4482E+00 | 5.0068E-06 | 3.8671E-04 | 4.6234E-02 | 2.5369E+01 | 4.1842E-05 | 2.7359E-04 |
| Proc 8 | 2.4557E-05 | 1.5258E-03 | 1.4909E+00 | 4.8876E-06 | 3.6669E-04 | 4.6220E-02 | 1.1133E+01 | 3.9220E-05 | 2.9099E-04 |
| Proc 9 | 6.1035E-05 | 1.9531E-03 | 5.3333E-01 | 7.6294E-06 | 4.9897E-04 | 8.5969E-02 | 5.0000E+00 | 9.1553E-05 | 2.1973E-03 |
| Proc 10 | 6.1035E-05 | 1.9531E-03 | 1.6523E+04 | 1.5259E-05 | 9.7656E-04 | 7.5120E-02 | 7.0000E+00 | 6.8665E-05 | 1.1902E-03 |
| Prov 11 | 6.1035E-05 | 1.9531E-03 | 1.6522E+04 | 7.6294E-06 | 9.7656E-04 | 7.5120E-02 | 7.0000E+00 | 7.6294E-05 | 1.1902E-03 |
| Proc 12 | 1.5259E-05 | 1.9531E-03 | 1.6522E+04 | 7.6294E-06 | 4.8828E-04 | 7.1746E-02 | 9.0000E+00 | 7.2241E-05 | 1.2207E-03 |

## Table 4.4: MAXARE under Double Precision

| MAXARE(D) | TP1 | TP2 | TP3 | TP41 | TP42 | TP43 | TP44 | TP5 | TP6 |
|---|---|---|---|---|---|---|---|---|---|
| Proc 1 | 1.1369E-13 | 2.4714E-12 | 1.8626E-09 | 1.4211E-14 | 1.8190E-12 | 1.1642E-10 | 1.4901E-08 | 1.1369E-13 | 5.1914E-13 |
| Proc 2 | 5.6843E-14 | 3.6380E-12 | 1.4461E-09 | 1.4211E-14 | 9.0949E-13 | 1.7462E-10 | 7.4506E-09 | 1.2790E-13 | 9.6634E-13 |
| Proc 3 | 2.2737E-13 | 5.4570E-12 | 2.5251E-09 | 6.1950E-14 | 6.3665E-12 | 3.4925E-10 | 3.4925E-10 | 2.7001E-13 | 4.2064E-12 |
| Proc 4 | 1.1369E-13 | 2.3647E-11 | 4.2285E-05 | 2.1316E-14 | 2.3099E-12 | 1.7462E-10 | 2.2352E-08 | 1.5632E-13 | 1.3642E-12 |
| Proc 5 | 2.8422E-13 | 4.7893E-12 | 1.9447E-09 | 1.5987E-14 | 2.3055E-12 | 2.3283E-10 | 2.9802E-08 | 8.8107E-13 | 1.6598E-11 |
| Proc 6 | 1.1369E-13 | 3.6380E-12 | 1.8626E-09 | 1.4211E-14 | 9.0949E-13 | 1.1642E-10 | 7.5181E-09 | 1.1369E-13 | 1.1369E-12 |
| Proc 7 | 1.1369E-13 | 3.6380E-12 | 1.8626E-09 | 1.4211E-14 | 1.8190E-12 | 1.7462E-10 | 1.4901E-08 | 1.7053E-13 | 2.1600E-12 |
| Proc 8 | 1.1369E-13 | 4.2912E-12 | 1.8626E-09 | 2.8422E-14 | 3.6380E-12 | 3.4925E-10 | 2.0740E-08 | 2.2737E-13 | 2.6148E-12 |
| Proc 9 | 1.1369E-13 | 3.6380E-12 | 1.4461E-09 | 1.4211E-14 | 1.8190E-12 | 1.1642E-10 | 7.4506E-09 | 1.5632E-13 | 3.6113E-12 |
| Proc 10 | 5.6843E-14 | 3.7313E-12 | 2.0940E-05 | 1.4211E-14 | 1.8190E-12 | 1.7462E-10 | 2.2352E-08 | 1.5632E-13 | 1.2186E-12 |
| Proc 11 | 5.6843E-14 | 3.7313E-12 | 2.0941E-05 | 1.4211E-14 | 1.8190E-12 | 1.7462E-10 | 2.2352E-08 | 1.7053E-13 | 1.2186E-12 |
| Proc 12 | 1.1369E-13 | 3.7313E-12 | 2.0941E-05 | 1.4211E-14 | 1.8190E-12 | 1.1642E-10 | 1.4901E-08 | 1.5632E-13 | 1.2186E-12 |

## Table 4.5: ORE under Single Precision

| ORE(S) | TP1 | TP2 | TP3 | TP41 | TP42 | TP43 | TP44 | TP5 | TP6 |
|---|---|---|---|---|---|---|---|---|---|
| Proc 1 | 7.9721E-05 | 2.8165E-02 | 3.7670E+00 | 2.1267E-04 | 1.7017E-02 | 1.8252E+00 | 5.2130E+02 | 7.7588E-02 | 1.0781E+01 |
| Proc 2 | 8.4311E-05 | 1.9942E-02 | 5.1353E+00 | 3.4666E-04 | 1.4591E-02 | 1.4978E+00 | 1.4851E+02 | 7.9670E-02 | 1.0825E+01 |
| Proc 3 | 2.6092E-04 | 3.5845E-02 | 5.4519E+00 | 2.6882E-04 | 2.8050E-02 | 1.6858E+00 | 1.6858E+00 | 7.8870E-02 | 1.1015E+01 |
| Proc 4 | 1.3366E-04 | 1.8962E-02 | 3.6441E+00 | 1.7709E-04 | 1.5828E-02 | 1.4867E+00 | NaN | 7.6250E-02 | 1.0542E+01 |
| Proc 5 | 1.7825E-04 | 1.8865E-02 | 5.8119E+00 | 1.4186E-04 | 9.3399E-03 | 1.1775E+00 | 3.3758E+02 | 7.5654E-02 | 1.0529E+01 |
| Proc 6 | 1.2526E-04 | 1.9844E-02 | 6.0251E+00 | 1.2279E-04 | 9.7828E-03 | 8.9603E-01 | 1.9242E+02 | 7.6298E-02 | 1.0542E+01 |
| Proc 7 | 2.4018E-04 | 1.9384E-02 | 5.7981E+00 | 1.3071E-04 | 1.0204E-02 | 9.7736E-01 | 2.3153E+02 | 7.6570E-02 | 1.0531E+01 |
| Proc 8 | 1.2872E-04 | 1.6985E-02 | 6.0274E+00 | 1.2934E-04 | 9.7471E-03 | 9.3048E-01 | 1.9473E+02 | 7.5713E-02 | 1.0526E+01 |
| Proc 9 | 2.0275E-04 | 1.6383E-02 | 3.5601E+00 | 1.1259E-04 | 5.3112E-03 | 7.8108E-01 | 8.5883E+01 | 1.3900E-01 | 3.1917E+00 |
| Proc 10 | 1.8367E-04 | 2.7213E-02 | 3.3883E+04 | 2.2340E-04 | 1.7498E-02 | 2.1203E+00 | 1.7003E+02 | 8.0323E-02 | 1.0972E+01 |
| Prov 11 | 1.7092E-04 | 2.7842E-02 | 3.3883E+04 | 2.4676E-04 | 1.7498E-02 | 2.1203E+00 | 1.7203E+02 | 7.9933E-02 | 1.0988E+01 |
| Proc 12 | 7.1913E-05 | 2.6956E-02 | 3.3883E+04 | 1.9002E-04 | 1.4920E-02 | 1.6139E+00 | 1.5876E+02 | 8.0110E-02 | 1.0927E+01 |

## Table 4.6: ORE under Double Precision

| ORE(D) | TP1 | TP2 | TP3 | TP41 | TP42 | TP43 | TP44 | TP5 | TP6 |
|---|---|---|---|---|---|---|---|---|---|
| **Proc 1** | 3.1486E-13 | 4.2940E-11 | 8.2400E-09 | 4.3465E-13 | 4.1593E-11 | 3.4351E-09 | 3.2901E-07 | 1.4603E-10 | 2.0170E-08 |
| **Proc 2** | 2.9809E-13 | 3.8647E-11 | 7.2145E-09 | 4.4076E-13 | 2.8369E-11 | 4.0932E-09 | 1.9688E-07 | 1.4910E-10 | 2.0279E-08 |
| **Proc 3** | 1.1076E-12 | 1.0909E-10 | 1.6599E-08 | 1.6265E-12 | 1.0069E-10 | 8.7043E-09 | 8.7043E-09 | 4.3246E-10 | 1.1646E-07 |
| **Proc 4** | 4.0834E-13 | 3.5271E-10 | 1.6354E-04 | 6.8556E-13 | 6.8057E-11 | 4.7508E-09 | 5.2569E-08 | 1.5393E-10 | 2.0470E-08 |
| **Proc 5** | 5.8542E-13 | 8.0202E-11 | 1.4783E-08 | 7.6406E-13 | 6.9908E-11 | 6.9648E-09 | 7.2536E-07 | 1.7831E-10 | 2.2359E-08 |
| **Proc 6** | 3.1319E-13 | 3.9636E-11 | 6.9062E-09 | 3.8036E-13 | 2.9394E-11 | 3.2306E-09 | 2.2832E-07 | 1.5009E-10 | 2.0217E-08 |
| **Proc 7** | 3.6315E-13 | 4.8937E-11 | 1.1378E-08 | 3.6848E-13 | 4.4544E-11 | 4.6604E-09 | 4.2368E-07 | 1.5067E-10 | 2.0382E-08 |
| **Proc 8** | 4.5819E-13 | 1.0442E-10 | 1.5714E-08 | 8.0280E-13 | 8.5687E-11 | 1.1774E-08 | 6.2413E-07 | 3.8972E-10 | 1.1759E-07 |
| **Proc 9** | 2.9554E-13 | 2.8481E-11 | 5.2755E-09 | 2.8832E-13 | 1.9558E-11 | 1.5769E-09 | 1.4170E-07 | 2.5795E-10 | 5.9011E-08 |
| **Proc 10** | 2.4492E-13 | 4.1265E-11 | 5.2111E-06 | 3.4472E-13 | 3.7945E-11 | 3.1286E-09 | 3.4516E-07 | 1.5399E-10 | 2.0373E-08 |
| **Proc 11** | 2.1694E-13 | 4.2175E-11 | 5.2112E-05 | 3.2230E-13 | 3.7490E-11 | 3.3033E-09 | 3.4516E-07 | 1.5455E-10 | 2.0385E-08 |
| **Proc 12** | 2.8488E-13 | 3.6427E-11 | 5.2112E-05 | 2.9698E-13 | 3.0710E-11 | 3.0476E-09 | 2.5663E-07 | 1.5188E-10 | 2.0271E-08 |

## Table 4.7 for MINAE(S, D) between Single and Double Precision calculations

| MINAE(SD) | TP1 | TP2 | TP3 | TP41 | TP42 | TP43 | TP44 | TP5 | TP6 |
|---|---|---|---|---|---|---|---|---|---|
| Proc 1 | 2.8422E-14 | 1.8878E-08 | 8.9733E-08 | 3.7478E-08 | 2.5415E-07 | 2.5415E-07 | 2.5415E-07 | 3.4895E-10 | 3.3754E-10 |
| Proc 2 | 4.4409E-16 | 3.7599E-05 | 2.9789E-02 | 1.0393E-07 | 8.1511E-09 | 3.4234E-07 | 2.18O4E-07 | 5.7804E-10 | 1.9259E-10 |
| Proc 3 | 4.6971E-08 | 3.8719E-05 | 1.4642E-02 | 4.9425E-08 | 1.1721E-04 | 1.4576E-03 | 1.4576E-03 | 4.2617E-09 | 4.8774E-09 |
| Proc 4 | 1.0398E-07 | 4.2477E-05 | 2.1705E-02 | 1.3248E-07 | 2.8595E-05 | 4.7612E-04 | Inf | 5.2125E-10 | 4.2748E-07 |
| Proc 5 | 9.5415E-09 | 3.0867E-04 | 3.0409E-03 | 2.5835E-08 | 3.8400E-06 | 4.0875E-04 | 7.5673E-02 | 5.1557E-11 | 1.3074E-12 |
| Proc 6 | 3.8137E-09 | 1.8545E-04 | 3.0409E-03 | 4.4225E-09 | 2.8390E-06 | 2.8088E-04 | 5.6237E-02 | 2.0264E-06 | 1.9400E-05 |
| Proc 7 | 3.0546E-08 | 3.3281E-04 | 3.0409E-03 | 3.8417E-07 | 5.2825E-05 | 5.2960E-03 | 6.7380E-01 | 1.6367E-07 | 1.8749E-07 |
| Proc 8 | 1.2641E-09 | 1.8556E-04 | 3.0409E-03 | 1.6998E-08 | 2.7935E-06 | 2.8083E-04 | 5.6237E-02 | 6.4659E-12 | 5.5138E-12 |
| Proc 9 | 0.0000E+00 | 7.8481E-08 | 1.1097E-07 | 1.0259E-08 | 1.5811E-08 | 2.3757E-08 | 5.6944E-09 | 2.2329E-09 | 2.5949E-09 |
| Proc 10 | 0.0000E+00 | 1.8878E-08 | 8.9733E-08 | 1.2903E-07 | 9.2766E-08 | 2.5415E-07 | 2.5415E-07 | 9.4606E-10 | 8.7670E-10 |
| Prov 11 | 0.0000E+00 | 1.8878E-08 | 8.9733E-08 | 1.2903E-07 | 9.2766E-08 | 2.5415E-07 | 2.5415E-07 | 9.4606E-10 | 1.4912E-09 |
| Proc 12 | 0.0000E+00 | 1.8878E-08 | 8.9733E-08 | 3.7478E-08 | 3.7478E-08 | 3.7478E-08 | 3.7478E-08 | 9.4609E-10 | 8.2019E-10 |

## Table 4.8: MAXAE(S, D) errors between Single and Double Precision calculations

| MAXAE(SD) | TP1 | TP2 | TP3 | TP41 | TP42 | TP43 | TP44 | TP5 | TP6 |
|---|---|---|---|---|---|---|---|---|---|
| Proc 1 | 2.7127E-05 | 4.5160E-01 | 1.3403E+05 | 5.8003E-05 | 6.8901E-02 | 6.7987E+02 | 1.0528E+08 | 8.0550E-05 | 1.2098E-03 |
| Proc 2 | 4.9506E-05 | 3.7231E+00 | 9.7505E+04 | 5.0969E-05 | 3.8604E-01 | 3.6333E+03 | 2.5358E+07 | 7.4475E-04 | 2.4924E-02 |
| Proc 3 | 6.9505E-03 | 3.2695E+00 | 5.5296E+05 | 1.4576E-04 | 1.0315E+00 | 7.0150E+03 | 7.0150E+03 | 2.1435E-05 | 3.8256E-03 |
| Proc 4 | 2.9093E-03 | 1.7307E+00 | 2.6785E+05 | 1.0287E-04 | 2.3300E+00 | 1.2735E+04 | Inf | 1.0047E-03 | 1.1509E-02 |
| Proc 5 | 2.8536E-03 | 1.8751E+00 | 1.8137E+05 | 5.3411E-05 | 3.4919E-01 | 3.6567E+03 | 1.2861E+08 | 1.3310E-06 | 4.3696E-06 |
| Proc 6 | 2.5588E-03 | 8.9710E-01 | 1.8137E+05 | 3.8305E-05 | 2.4802E-01 | 2.6363E+03 | 4.9001E+07 | 4.5307E-06 | 2.9365E-05 |
| Proc 7 | 2.2191E-03 | 2.0862E+00 | 1.8137E+05 | 5.7456E-05 | 3.4977E-01 | 3.6567E+03 | 1.2861E+08 | 3.2308E-04 | 1.2438E-02 |
| Proc 8 | 2.5609E-03 | 8.9731E-01 | 1.8137E+05 | 3.9159E-05 | 2.4806E-01 | 2.6363E+03 | 4.9001E+07 | 8.7661E-07 | 2.6700E-06 |
| Proc 9 | 6.7817E-03 | 2.2032E-03 | 1.2214E+00 | 1.4700E-05 | 1.0511E-03 | 1.1536E-01 | 7.6999E+06 | 8.5831E-05 | 2.2106E-03 |
| Proc 10 | 8.8162E-05 | 1.4282E+00 | 1.8206E+05 | 1.9367E-05 | 6.0751E-02 | 6.1944E+02 | 7.3597E+06 | 1.3911E-04 | 2.4661E-03 |
| Proc 11 | 8.8162E-05 | 1.4282E+00 | 1.8206E+05 | 2.0312E-05 | 6.0751E-02 | 6.1944E+02 | 7.3597E+06 | 1.3911E-04 | 2.4356E-03 |
| Proc 12 | 2.7127E-05 | 1.4272E+00 | 1.8206E+05 | 1.7415E-05 | 6.0359E-02 | 6.1942E+02 | 7.3597E+06 | 1.2385E-04 | 2.4356E-03 |

**Table 4.9: REL(S, D) errors between Single and Double Precision calculations**

| REL(S, D) | TP1 | TP2 | TP3 | TP41 | TP42 | TP43 | TP44 | TP5 | TP6 |
|---|---|---|---|---|---|---|---|---|---|
| Proc 1 | 1.8238E-04 | 9.7760E+00 | 8.5381E+05 | 2.0605E-03 | 2.9959E+00 | 2.9939E+04 | 4.6364E+09 | 1.8035E-01 | 3.7116E+01 |
| Proc 2 | 2.2072E-04 | 5.5973E+01 | 6.8745E+05 | 9.9624E-04 | 9.6408E+00 | 9.0833E+04 | 6.3394E+08 | 2.5728E-01 | 3.8825E+01 |
| Proc 3 | 4.1639E-02 | 5.6868E+01 | 3.7384E+06 | 4.9278E-03 | 3.4096E+01 | 2.2813E+05 | 2.2813E+05 | 3.9122E-01 | 1.7572E+02 |
| Proc 4 | 1.7467E-02 | 1.8265E+01 | 1.5096E+06 | 2.2607E-03 | 3.1956E+01 | 2.1517E+05 | NaN | 2.3411E-01 | 3.2271E+01 |
| Proc 5 | 1.7052E-02 | 2.9542E+01 | 1.5219E+06 | 1.7384E-03 | 1.2950E+01 | 1.3674E+05 | 3.9627E+09 | 2.5296E-03 | 2.0951E-01 |
| Proc 6 | 1.5258E-02 | 1.8216E+01 | 1.5219E+06 | 1.4254E-03 | 1.0783E+01 | 1.1486E+05 | 2.1067E+09 | 2.9560E-02 | 5.9796E+00 |
| Proc 7 | 1.3236E-02 | 3.3603E+01 | 1.5219E+06 | 1.5823E-03 | 1.0819E+01 | 1.1517E+05 | 3.4550E+09 | 8.9744E-02 | 1.8154E+01 |
| Proc 8 | 1.5273E-02 | 1.8219E+01 | 1.5219E+06 | 1.4598E-03 | 1.0785E+01 | 1.1486E+05 | 2.1067E+09 | 2.2196E-03 | 1.4836E-01 |
| Proc 9 | 2.9606E-04 | 2.6209E-02 | 5.3162E+00 | 2.4325E-04 | 1.3152E-02 | 1.1988E+00 | 1.2451E+02 | 1.3272E-01 | 3.1786E+01 |
| Proc 10 | 4.8392E-04 | 4.0830E+01 | 1.5125E+06 | 4.9828E-04 | 2.6174E+00 | 2.7278E+04 | 3.2413E+08 | 2.3295E-01 | 6.1311E+01 |
| Proc 11 | 4.9096E-04 | 4.0828E+01 | 1.5125E+06 | 4.7100E-04 | 2.6174E+00 | 2.7278E+04 | 3.2413E+08 | 2.3367E-01 | 6.1499E+01 |
| Proc 12 | 1.6090E-04 | 4.0822E+01 | 1.5125E+06 | 4.4613E-04 | 2.6199E+00 | 2.7278E+04 | 3.2413E+08 | 1.9630E-01 | 5.7349E+01 |

**Table 4.10: Average extra digits between Single and Double Precision calculations**

| ANED | TP1 | TP2 | TP3 | TP41 | TP42 | TP43 | TP44 | TP5 | TP6 |
|---|---|---|---|---|---|---|---|---|---|
| Proc 1 | *** 6.7522 | 4.9691 | 3.2245 | 6.5086 | 5.0423 | 3.2272 | 0.7819 | 6.8841 | 6.7137 |
| Proc 2 | ** 7.1182 | 4.5133 | 2.4592 | 6.7433 | 5.0996 | 3.5632 | 2.3108 | 6.9038 | 6.8656 |
| Proc 3 | 6.2773 | 3.9600 | 1.7416 | 6.0373 | 3.9882 | 2.1693 | 2.1693 | 6.6245 | 6.0247 |
| Proc 4 | 6.3811 | 4.5148 | 1.8947 | 6.3439 | 4.0975 | 2.1999 | NaN | 7.0449 | 6.8526 |
| Proc 5 | 6.8213 | 4.3463 | 2.1006 | 6.5652 | 4.5249 | 2.5096 | 0.1074 | 8.7662 | 8.9569 |
| Proc 6 | 6.7783 | 4.3526 | 2.1006 | 6.7684 | 4.8025 | 2.7873 | 0.5022 | 7.5321 | 7.3209 |
| Proc 7 | 6.8009 | 4.1487 | 2.1006 | 6.4926 | 4.4743 | 2.4530 | 0.1447 | 7.3531 | 7.2240 |
| Proc 8 | 6.9079 | 4.3525 | 2.1006 | 6.7396 | 4.8026 | 2.7873 | 0.5022 | 8.8273 | 9.0882 |
| Proc 9 | * 7.3504 | 7.2928 | 7.3526 | 7.3681 | 7.4157 | 7.4296 | 7.3321 | 7.0729 | 6.8100 |
| Proc 10 | * 6.6272 | 4.5093 | 2.9679 | 6.8748 | 5.1722 | 3.3115 | 1.4794 | 6.8250 | 6.4886 |
| Proc 11 | * 6.5949 | 4.5091 | 2.9679 | 6.9339 | 5.1722 | 3.3115 | 1.4794 | 6.8266 | 6.4853 |
| Proc 12 | *** 6.7660 | 4.5436 | 2.9679 | 7.0429 | 5.5830 | 4.1590 | 2.7397 | 6.9129 | 6.5384 |

**Table 4.11: PCZE(S)**

| PCZE(S) | TP1 | TP2 | TP3 | TP41 | TP42 | TP43 | TP44 | TP5 | TP6 |
|---|---|---|---|---|---|---|---|---|---|
| Proc 1 | 41.67% | 14.06% | 24.00% | 28.00% | 22.00% | 22.00% | 0.00% | 18.75% | 17.50% |
| Proc 2 | 58.33% | 21.88% | 16.00% | 18.00% | 24.00% | 25.00% | 3.00% | 18.81% | 17.49% |
| Proc 3 | 5.56% | 0.00% | 0.00% | 0.00% | 0.00% | 0.00% | 0.00% | 0.56% | 0.12% |
| Proc 4 | 13.89% | 0.00% | 0.00% | 2.00% | 0.00% | 0.00% | 0.00% | 0.58% | 0.12% |
| Proc 5 | 11.11% | 0.00% | 0.00% | 4.00% | 0.00% | 0.00% | 0.00% | 0.60% | 0.11% |
| Proc 6 | 22.22% | 0.00% | 0.00% | 8.00% | 0.00% | 0.00% | 0.00% | 0.57% | 0.12% |
| Proc 7 | 16.67% | 0.00% | 0.00% | 3.00% | 0.00% | 0.00% | 0.00% | 0.58% | 0.11% |
| Proc 8 | 16.67% | 0.00% | 0.00% | 5.00% | 0.00% | 0.00% | 0.00% | 0.75% | 0.12% |
| Proc 9 | 36.11% | 18.75% | 16.00% | 36.00% | 40.00% | 32.00% | 11.00% | 10.80% | 7.85% |
| Proc 10 | 47.22% | 20.31% | 16.00% | 25.00% | 25.00% | 15.00% | 0.00% | 18.83% | 17.51% |
| Proc 11 | 44.44% | 17.19% | 16.00% | 23.00% | 25.00% | 16.00% | 0.00% | 19.32% | 17.47% |
| Proc 12 | 41.67% | 12.50% | 16.00% | 33.00% | 28.00% | 30.00% | 7.00% | 18.63% | 17.47% |

**Table 4.12: PCZE(D)**

| PZE(D) | TP1 | TP2 | TP3 | TP41 | TP42 | TP43 | TP44 | TP5 | TP6 |
|---|---|---|---|---|---|---|---|---|---|
| Proc 1 | 47.22% | 9.38% | 20.00% | 24.00% | 19.00% | 17.00% | 25.00% | 18.86% | 17.45% |
| Proc 2 | 33.33% | 21.88% | 24.00% | 21.00% | 23.00% | 22.00% | 32.00% | 18.90% | 17.38% |
| Proc 3 | 5.56% | 7.81% | 12.00% | 7.00% | 7.00% | 12.00% | 9.00% | 7.10% | 3.45% |
| Proc 4 | 25.00% | 4.69% | 8.00% | 14.00% | 9.00% | 22.00% | 19.00% | 18.79% | 17.34% |
| Proc 5 | 36.11% | 7.81% | 12.00% | 10.00% | 17.00% | 13.00% | 18.00% | 18.20% | 16.94% |
| Proc 6 | 52.78% | 17.19% | 16.00% | 31.00% | 23.00% | 23.00% | 25.00% | 18.34% | 17.39% |
| Proc 7 | 38.89% | 14.06% | 8.00% | 28.00% | 22.00% | 17.00% | 19.00% | 18.64% | 17.51% |
| Proc 8 | 27.78% | 7.81% | 4.00% | 17.00% | 11.00% | 11.00% | 12.00% | 7.74% | 2.95% |
| Proc 9 | 47.22% | 23.44% | 28.00% | 34.00% | 27.00% | 31.00% | 34.00% | 11.85% | 7.07% |
| Proc 10 | 41.67% | 20.31% | 16.00% | 28.00% | 23.00% | 22.00% | 23.00% | 18.29% | 17.29% |
| Proc 11 | 44.44% | 20.31% | 16.00% | 32.00% | 24.00% | 20.00% | 23.00% | 18.01% | 17.40% |
| Proc 12 | 55.56% | 23.44% | 16.00% | 33.00% | 25.00% | 20.00% | 31.00% | 18.52% | 17.43% |

# Appendix 5: Charts of selected accuracy measures for all Test Problems

## Chart 5.1: MAXARE(S)

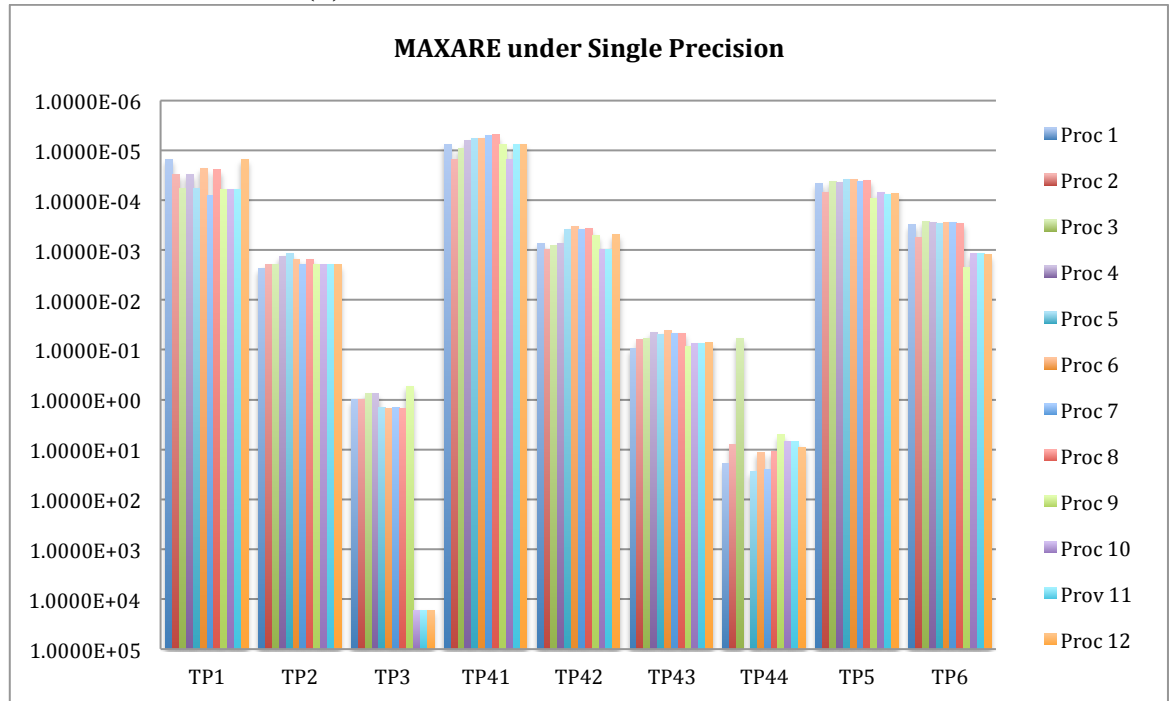

## Chart 5.2: MAXARE(D)

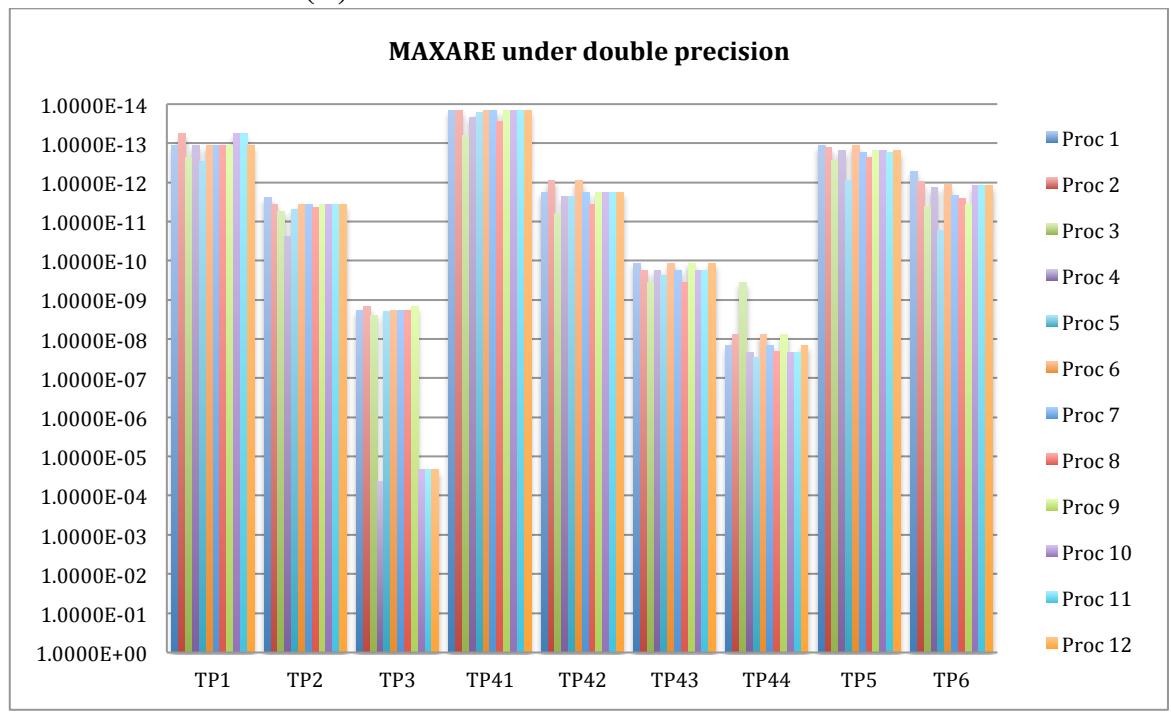

**Chart 5.3 : ORE(S)**

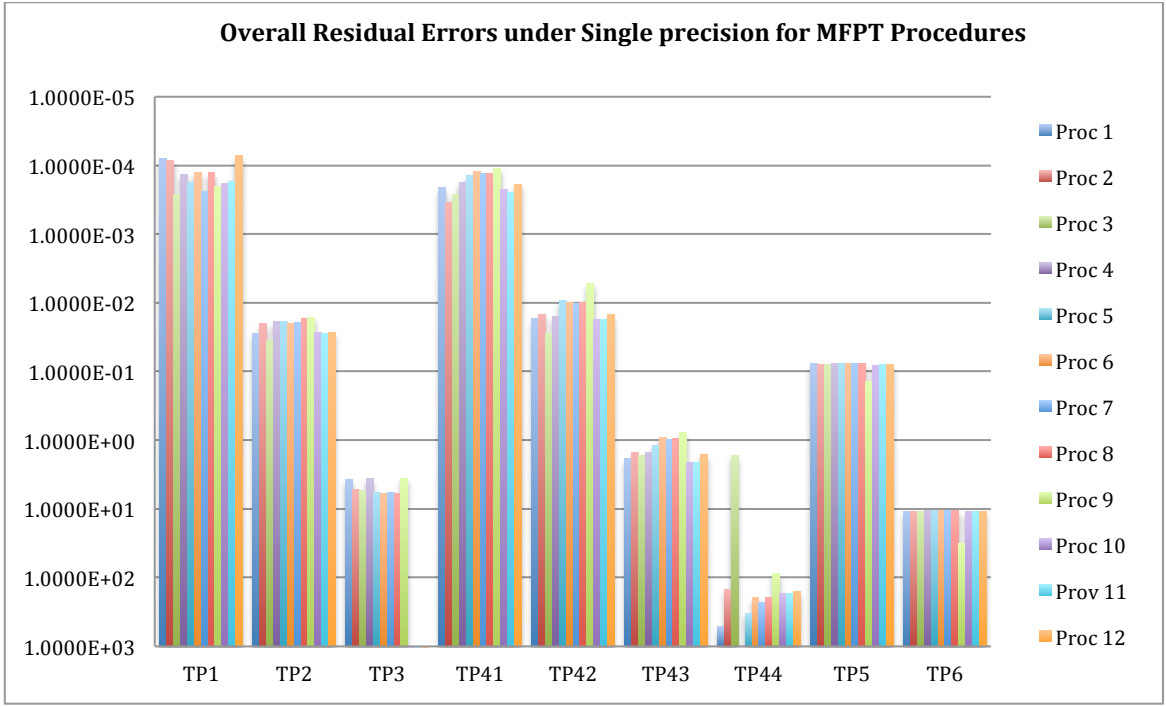

**Chart 5.4: ORE(D)**

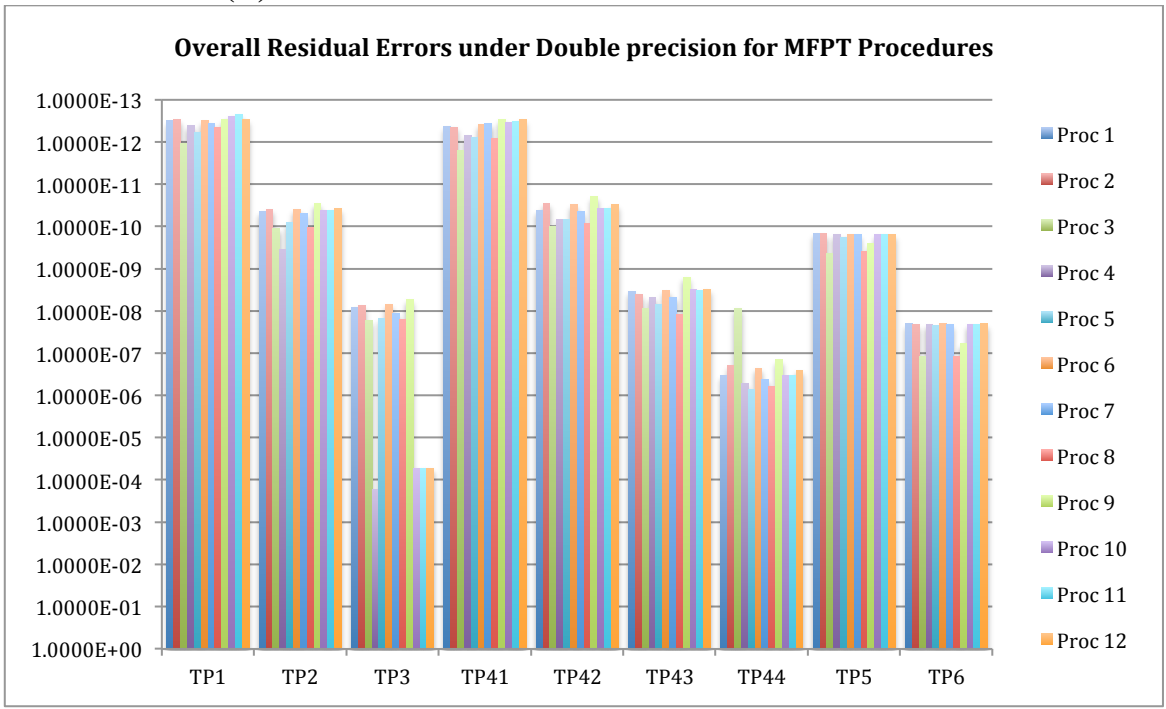

**Chart 5.5: REL(S,D)**

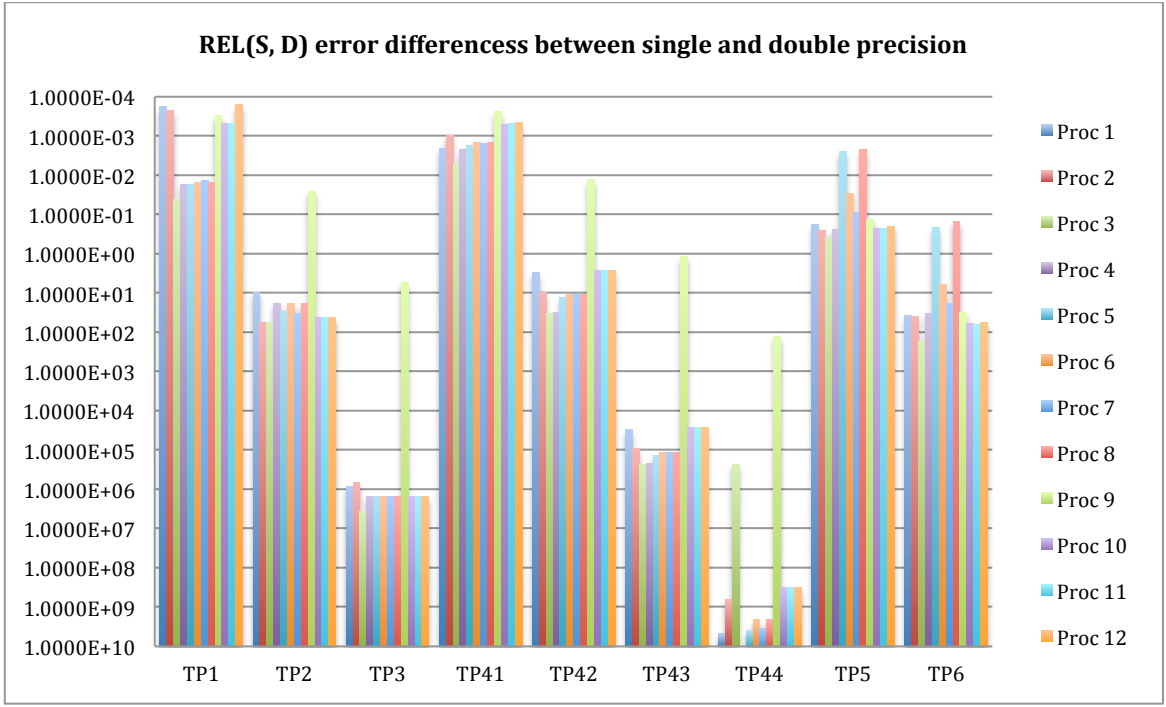

**Chart 5.6: ANED**

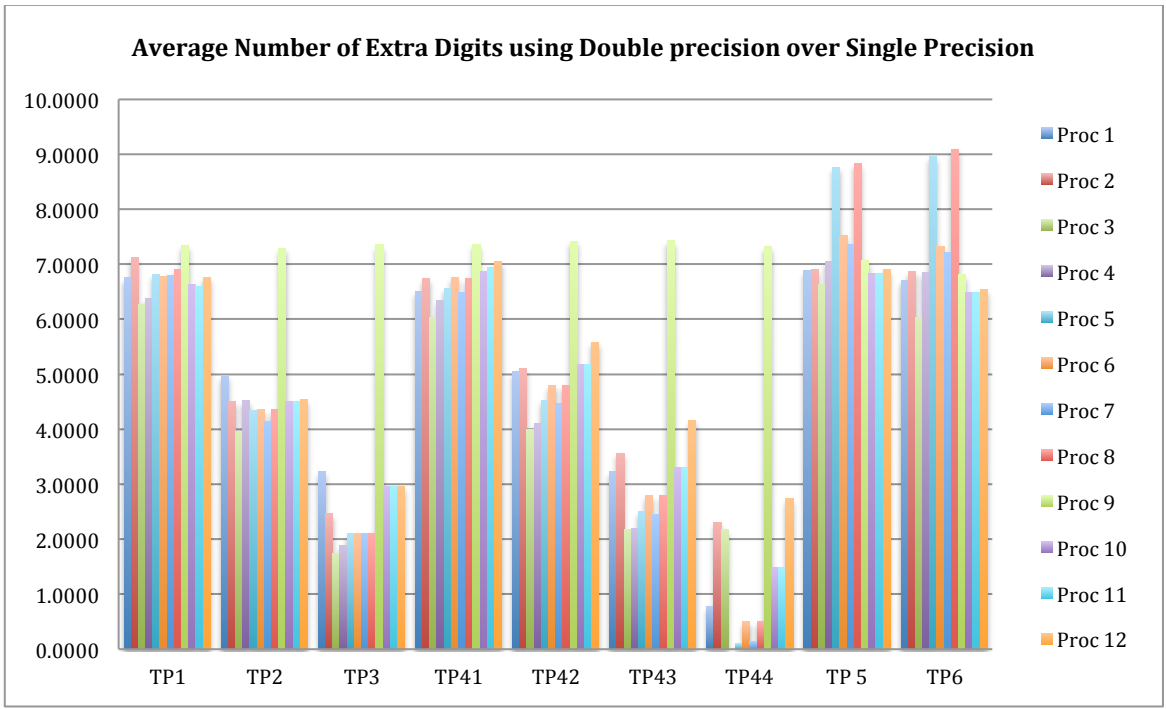

**Chart 5.7: PZCE(S)**

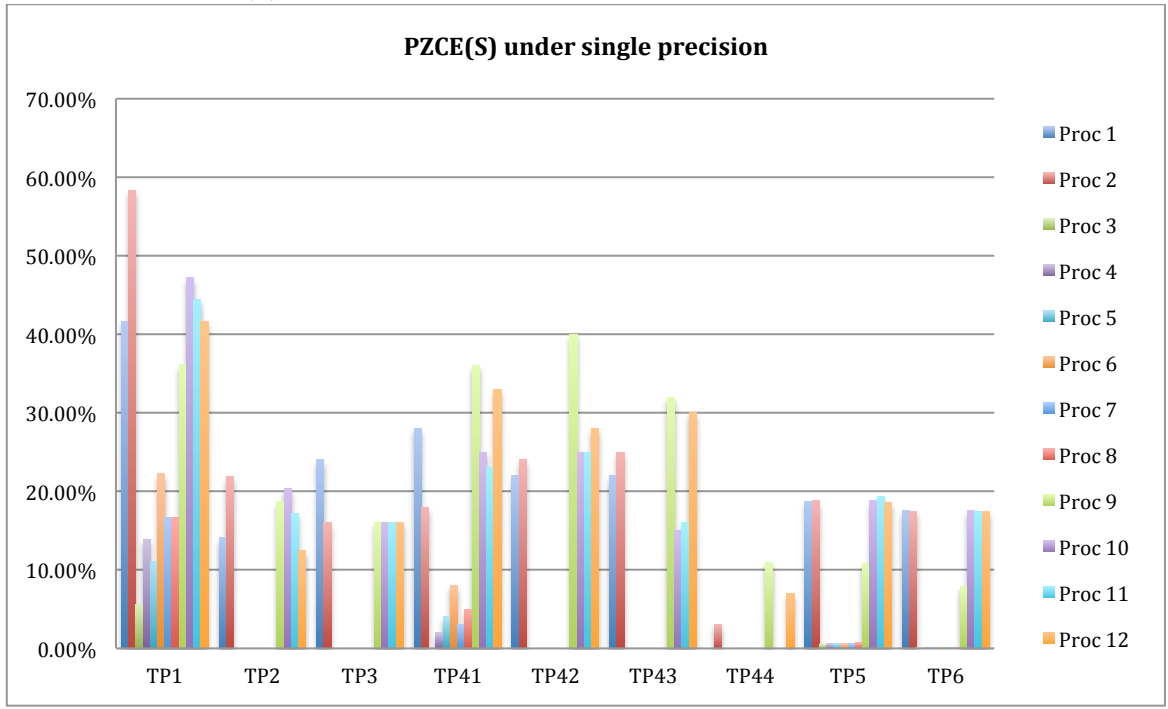

**Chart 5.8: PZCE(D)**

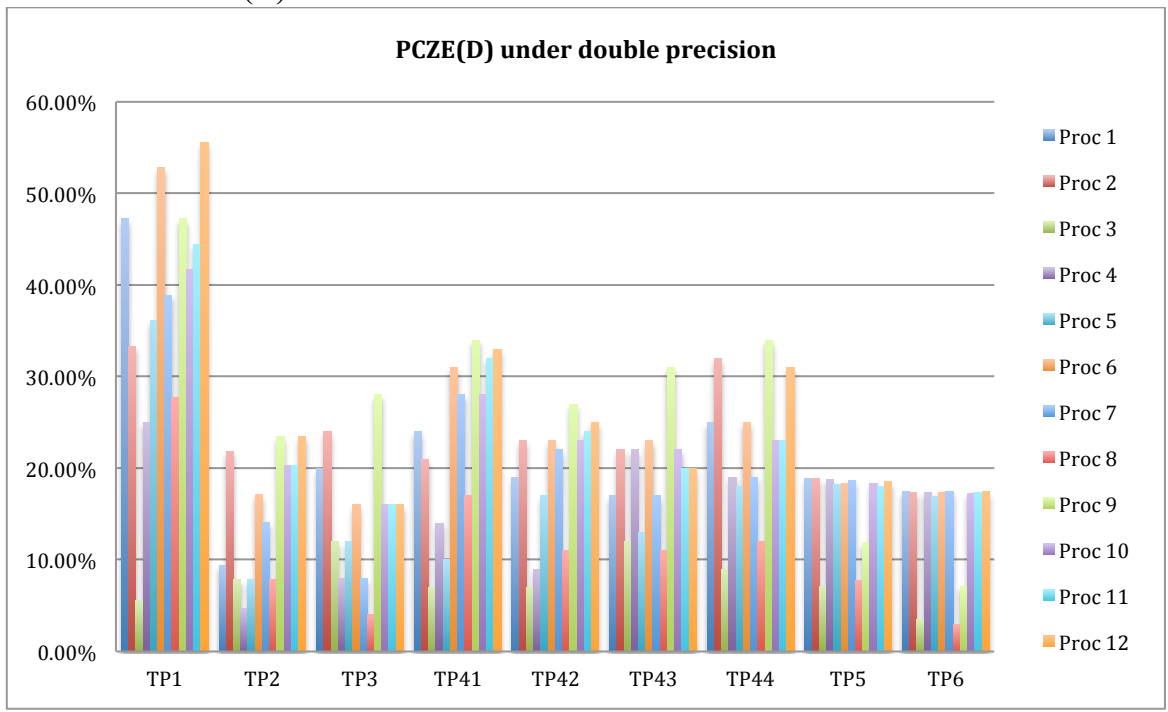